\documentclass[12pt]{report}
\usepackage{amsmath, amsthm, amssymb, epsf, graphicx}

\usepackage{uga}

\newtheorem*{thm}{Theorem}
\newtheorem*{lem}{Lemma}

\title{Random Hyperplanes Of A Convex Body, Sylvester's Problem And Crofton's Formula}

\author{Sheree T. Sharpe}

\professor{Malcolm R. Adams}

\fulltitle{Random Hyperplanes Of A Convex Body, Sylvester's Problem And Crofton's Formula}

\indexwords{Crofton's~Formula, Convex~Body}

\begin{document}
\pagenumbering{roman}

\begin{abstract}
Motivated by a problem on the 67th William Lowell Putnam Mathematical Competition,
we will summarize three different solutions found on a website. This Putman problem is a special case of Sylvester's four point problem! Suppose four points are taken at random in a convex body; what is the probability that they form a convex quadrilateral? We will see that there exists a relationship among Crofton's formula, random secants in two dimensions and the solution to this question. We will then present the solution following Kingman [3] to the Sylvester's four point problem in two and three dimensions for a unit ball by looking at convex bodies in three and four dimensions, respectively.
\end{abstract}

\tableofcontents

\listoffigures   % Optional - Omit this line if you don't want a list of figures.

\newpage
\pagenumbering{arabic}  % Ordinary pages have Arabic numerals.

\chapter{Introduction}

\noindent In the late nineteenth century, the English Mathematician
James Joseph Sylvester ($1814-1897$) proposed what is now considered a
classical problem of Geometrical Probability (which is, in more modern days,
referred to as Integral Geometry) in the Educational Times of 1864.
This problem is now known as the famous ``Sylvester Four Point
Problem," which states:
\begin{quote}
\textit{Find the probability that four points chosen at random
inside a convex set form a convex quadrilateral; that is, that none
of the points is inside the triangle formed by the other three
points.}
\end{quote}
The solution of Sylvester's four point problem uses techniques related to the proof of Crofton's formula. Crofton's formula is another classic result of Integral Geometry relating the length of a curve to the expected number of times a random line intersects the curve. Crofton's formula was named after the English Mathematician Morgan W. Crofton (1826-1915), one of the founders of Geometrical Probability. His methods allow the evaluation of some definite integrals without directly performing the integrations.

In this thesis, we will describe the $n$-dimensional results of John Frank Charles Kingman's 1969 paper ``Random Secants of a Convex Body" in Euclidean spaces of dimension $n=2$, $n=3$ and $n=4$. In this paper, Kingman computes the distribution of the $r$-dimensional affine subspace determined by randomly choosing $r+1 \leq n$ points in an $n$-dimensional convex body. For example, if we let two points be taken at random in a $2$-dimensional convex body then we can find the distribution of the line joining the two points and compare this distribution with other distributions for random secants in the $2$-dimensional convex body. The results find application to the $(n-1)$-dimensional case of the problem of Sylvester for a unit ball.

The motivation for this thesis is a question on the 2006 William Lowell Putnam Mathematical Competition which is a special case of Sylvester's problem. In the
next section we will present this problem and three solutions that
appeared on the website [2]. Next we will look at Crofton's formulas and
their proof, and finally we will explore distributions and invariant
measures. In chapter two, we will present Kingman's results in
dimensions two, three, and four and see how they solve Sylvester's
problem for the unit ball.

\section{Putnam Problem And Three Solutions}

\noindent This thesis was motivated by a question on the 2006
William Lowell Putnam Mathematical Competition for undergraduate
students which stated:
\begin{quote}
\textit{Four points are chosen uniformly and independently at random
in the interior of a given circle. Find the probability that they
are the vertices of a convex quadrilateral.}
\end{quote}
Three different solutions were found on the website [2]. All three solutions reduce the question to computing the average area of a triangle formed by three random points inside a unit circle. Since the problem is independent of scale, we assume, without loss of generality, that the circle has radius $1$. We may ignore the case where three of the points are collinear, as this occurs with probability zero. With this assumption, the only way the four points can fail to form the vertices of a convex quadrilateral is if one of the points lies inside the triangle formed by the other three points. There are four such configurations, depending on which point lies inside the triangle, and they are mutually exclusive, because two points cannot both lie inside the triangle at the same time. Thus, the desired probability is $1$ minus $4$ times the probability that a point S lies inside $\triangle PQR$ formed by the three points $P, Q, R$. The probability that S lies inside the $\triangle PQR$ is the expected value of the area of $\triangle PQR$, (denoted $[PQR]$), divided by the area of the circle. The three solutions found the expected area of the $\triangle PQR$ to be $\frac{35}{48\pi}$. The area of the unit circle is $\pi$. Therefore, the answer to the Putnam question is $1-\dfrac{35}{12\pi^2}$.

Before we explain the three different solutions, let's define some concepts. A random variable is a function that associates a unique numerical value with every outcome of an experiment. There are two types of random variables, discrete and continuous. The probability distribution of a discrete random variable is a list of probabilities associated with its possible values. A continuous random variable is one which takes an uncountably many values. The expected value (expectation) of a random variable indicates its average or central value. Each random variable has a cumulative distribution function. This function gives the probability that the random variable $X$ is less than or equal to $x$, for every value $x$. For a continuous random variable, the cumulative distribution function is the integral of its probability density function, i.e. the probability density function of a continuous random variable is a function which can be integrated to obtain the probability that the random variable takes a value in a given interval.

\subsection{The first solution}

\noindent The first solution, by Daniel Kane, is summarized below:

Let $O$ denote the center of the circle, and let $P',Q',R'$ be the projections of $P,Q,R$ onto the circle from $O$. Mr. Kane first notes that there are three cases to consider.

Case $i)$: If $P',Q',R'$ lie on no semicircle together,
then $O$ is inside the triangle $PQR$, so
$$
[PQR]=[OPR]+[OQR]+[OPQ]
$$

\begin{figure}[h]
  \begin{center}
  \includegraphics[width=2in]{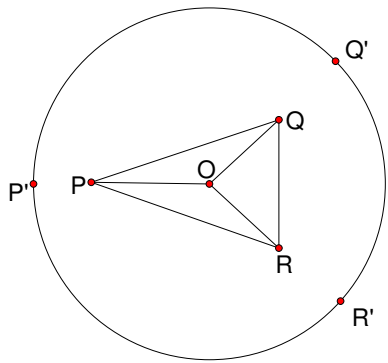}
  \caption{Putnam First Solution Case $i$}\label{Putnam First Solution Case $i$}
  \end{center}
\end{figure}

\pagebreak

Case $ii)$: If $P',Q',R'$ lie on the same semicircle in
that order and $Q$ lies inside $\triangle OPR$, then
$$
[PQR]=[OPR]-[OQR]-[OPQ]
$$

\begin{figure}[h]
  \begin{center}
  \includegraphics[width=2in]{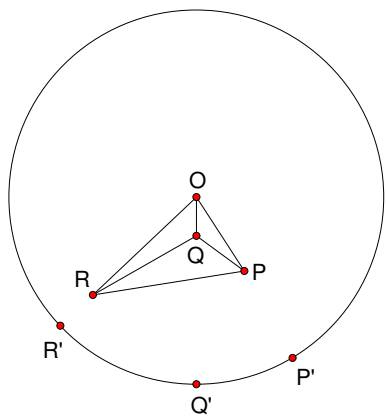}
  \caption{Putnam First Solution Case $ii$}\label{Putnam First Solution Case $ii$}
  \end{center}
\end{figure}

Case $iii)$: If $P',Q',R'$ lie on the same semicircle in
that order and $Q$ lies outside the triangle $OPR$, then
$$
[PQR]=-[OPR]+[OQR]+[OPQ]
$$

\begin{figure}[h]
  \begin{center}
  \includegraphics[width=2in]{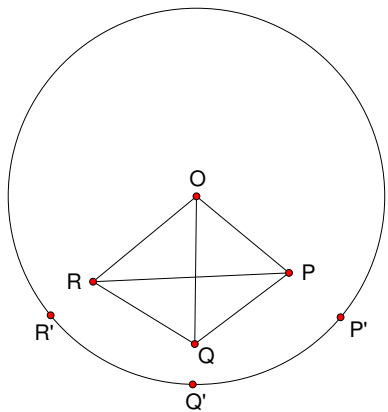}
  \caption{Putnam First Solution Case $iii$}\label{Putnam First Solution Case $iii$}
  \end{center}
\end{figure}

We first calculate the area $[OPQ]$ in all three cases.

\begin{figure}[h]
  \begin{center}
  \includegraphics[width=2in]{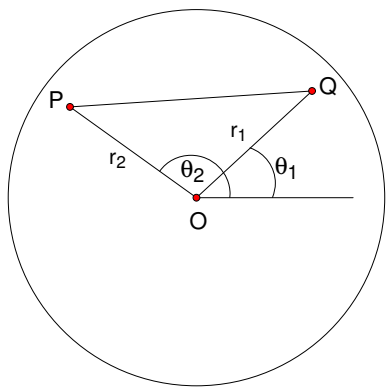}
  \caption{Area of triangle OPQ}\label{Area of triangle OPQ}
  \end{center}
\end{figure}

Write $r_1=OQ, r_2=OP$, and $\angle POQ = \theta_2 - \theta_1$. Then the area of triangle $OPQ$ is $[OPQ] = \frac{1}{2}r_1r_2\sin(\theta_2 - \theta_1)$. We can specify point $Q$ by $r_1$ and the angle $\theta_1$ it makes with the positive $x$-axis. We can specify point $P$ by $r_2$ and the angle $\theta_2$ it makes with the positive x-axis. Note that the distribution of $r_1$ is the area of the circle with radius $r_1$ divided by the area of the unit circle, i.e. $\frac{\pi r_1^2}{\pi}$. Thus the density of $r_1$ is given by $2r_1 dr_1$ on $[0,1]$, and similarly for $r_2$. The distribution of $\theta_1$ is uniform on $[0,2\pi]$ and is given by $\frac{1}{2\pi}d\theta_1$, and similarly for $\theta_2$. These distributions and densities are independent. Thus,
\begin{align*}
E([OPQ]) & = \int_0^{2\pi} \int_0^{2\pi} \int_0^1 \int_0^1 \frac{1}{2} r_1 r_2 |\sin(\theta_2 - \theta_1)| 2r_1 dr_1 2r_2 dr_2 \frac{1}{2\pi}d\theta_1 \frac{1}{2\pi}d\theta_2\\
& = \frac{1}{2\pi^2} \int_0^{2\pi} \int_0^{2\pi} \left( \int_0^1
r^2 dr\right)^2 |\sin(\theta_2 - \theta_1)| d\theta_1 d\theta_2\\
& = \frac{1}{18\pi^2} \int_0^{2\pi} \int_0^{2\pi}
|\sin(\theta_2 - \theta_1)| d\theta_1 d\theta_2.
\end{align*}
Holding $\theta_2$ fixed, let $x=\theta_1 - \theta_2$, so $dx= d\theta_1$. Since $0 \leq \theta_1 \leq 2\pi$, we have $-\theta_2 \leq x \leq 2\pi - \theta_2$. Then,
$$
\int_0^{2\pi} |\sin(\theta_2 - \theta_1)| d\theta_1 = \int_{-\theta_2}^{2\pi - \theta_2}|\sin(x)| dx = \int_0^{2\pi}|\sin(x)| dx =4
$$
Hence 
$$
E([OPQ]) = \frac{1}{18\pi^2} \int_0^{2\pi} 4 d\theta_2 = \frac{4}{9\pi}.
$$

We can now calculate the expected value, $E([PQR])$, of the area
$[PQR]$ in Case $i)$. By symmetry, we note that
$$
E([PQR]) = E([OPR]+[OQR]+[OPQ])=3E([OPQ])= \frac{4}{3\pi}.
$$

Cases $ii)$ and $iii)$ require a bit more care since the symmetry is broken.

Write $E '(X)$ for the expectation of a random variable $X$ restricted to case $ii)$ and case $iii)$. Let $\chi$ be the random variable with value $1$ if $Q$ is inside the $\triangle OPR$ and $0$ otherwise.
So the expected area of the $\triangle OPR$ given that $Q$ is inside the $\triangle OPR$ is
$$
E'([OPR]\cdot \chi) = E'\left(\dfrac{2[OPR]^2}{\theta}\right)
$$
The random variable, $\chi$, has value $1$ if $Q$ is inside the $\triangle OPR$, therefore the probability that $Q$ is inside the triangle is $\dfrac{[OPR]}{{\theta}/2}$, where ${\theta}/2$ is the area of the sector $OP'R'$.

\pagebreak

\begin{figure}[h]
  \begin{center}
  \includegraphics[width=2in]{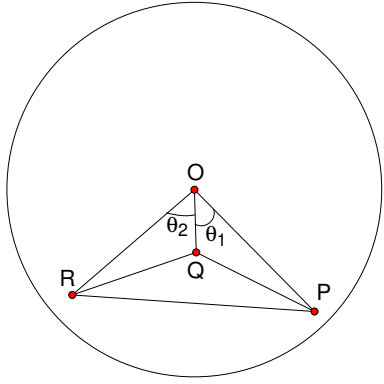}
  \caption{Q inside triangle OPR}\label{Q inside triangle OPR}
  \end{center}
\end{figure}

Let $\theta_1 = \angle POQ$ and $\theta_2 = \angle QOR$; then the distribution of $\theta_1, \theta_2$ is uniform on the region $0 \leq \theta_1$, $0 \leq \theta_2$, $\theta_1 + \theta_2 \leq \pi$. To compute the distribution for $\theta=\theta_1 + \theta_2$, note that the area of the triangle with base and height $s$ is $\frac{s^2}{2}$, and the area of the triangle with base and height $\pi$ is $\frac{\pi^2}{2}$. Thus, the probability that $s \leq \pi$ is $\frac{\frac{s^2}{2}}{\frac{\pi^2}{2}}= \frac{s^2}{\pi^2}$. In particular, the density of $\theta(s)$ is $\frac{2s}{\pi^2}$ on $[0,\pi]$. By the standard abuse of notation this density is written $\frac{2\theta}{\pi^2}d\theta$.
\begin{align*}
E'([OPR]\cdot \chi) & =E'\left(\dfrac{2\left(\frac{1}{2}r_Pr_R\sin(\theta)\right)^2}{\theta}\right)\\
 & = E'\left(\dfrac{1}{2}r_P^2r_R^2 \theta^{-1}\sin^2(\theta)\right)\\
 & =\frac{1}{2}\left(\int_0^1 2r r^2 dr\right)^2 \left(\int_0^{\pi} \dfrac{2\theta}{\pi^2}\theta^{-1}\sin^2(\theta)d\theta\right)\\
 & =\dfrac{1}{2}\left(\int_0^1 2r^3 dr\right)^2 \left(\int_0^{\pi} \frac{2}{\pi^2}\sin^2(\theta)d\theta\right)\\
 & =\dfrac{1}{8\pi}.
\end{align*}
Given $\triangle ABC$, if a point $D$ is chosen uniformly at random
inside $\triangle ABC$, then the expectation of $[DAB]$ is the
area of the triangle bounded by $AB$ and the centroid of $ABC$,
given by $\frac{1}{3}[ABC]$. 
So now in Case $ii)$,
\begin{align*}
E'([PQR]\cdot \chi)& =E'([OPR]\cdot \chi)-E'(([OQR]+[OQP])\cdot \chi)\\
& =E'([OPR]\cdot \chi)-\frac{2}{3}E'([OPR]\cdot \chi)\\
& =\frac{1}{3} \left(\frac{1}{8\pi}\right)\\
& =\frac{1}{24\pi}.
\end{align*}
And in Case $iii)$,
\begin{align*}
E'((1-\chi)[PQR])& =E'((1-\chi)([OPQ]+[OQR]-[OPR]))\\
& =E'([OPQ])+E'([OQR])-E'([OPR])\\
& -E'(([OPQ]+[OQR])\cdot \chi)+ E'([OPR]\cdot \chi)\\
& =E'([OPR])- \frac{2}{3} E'([OPR]\cdot \chi) + E'([OPR]\cdot \chi)\\
& =\frac{4}{9\pi} - \frac{2}{3} \left(\frac{1}{8\pi}\right) + \frac{1}{8\pi}\\
& =\frac{4}{9\pi} + \frac{1}{24\pi}.
\end{align*}
\noindent Then summing case $ii)$ and case $iii)$ gives $\frac{1}{24\pi} +
\frac{4}{9\pi} + \frac{1}{24\pi} = \frac{19}{36\pi}$.

The probability that $P', Q', R'$ lie on the same semicircle proceeding clockwise from $P'$ is $(\frac{1}{2})^2 = \frac{1}{4}$. Finally, $P', Q', R'$ lie on the same semicircle in some order occurs with probability $\frac{3}{4}$.
Hence case $i)$ occurs with probability $\frac{1}{4}$, and cases $ii)$ and $iii)$ occur together with probability $\frac{3}{4}$. That is,
$$
E([PQR]) = \frac{1}{4} \left(\frac{4}{3\pi}\right) + \frac{3}{4}\left(\frac{19}{36\pi}\right) = \frac{35}{48\pi}.
$$

\subsection{The second solution}

\noindent The second solution, by David Savitt, is summarized below:

Draw the lines $PQ,QR, RP$, which with probability $1$ divide the interior of the circle into seven regions. Put $a=[PQR]$, and let $b_1,b_2,b_3$ denote the areas of the three regions sharing a side with the $\triangle PQR$. Let $c_1,c_2,c_3$ denote the areas of the other three regions.

\begin{figure}[h]
  \begin{center}
  \includegraphics[width=2in]{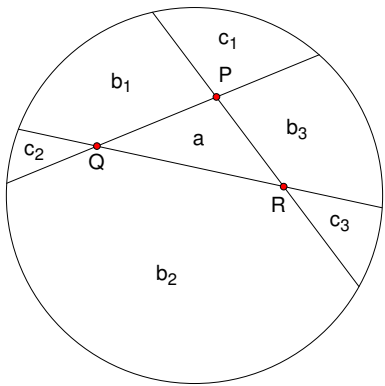}
  \caption{Putnam Second Solution}\label{Putnam Second Solution}
  \end{center}
\end{figure}

Put $A=E(a),B=E(b_i),C=E(c_i)$ where $i=1,2,3$, so that $A+3B+3C = \pi$, which is the area of the circle. Now, $c_1+c_2+c_3+a$ is the area of the region in which we can choose a fourth point $S$ so that the quadrilateral $PQRS$ fails to be convex. So the probability that the four points fail to form a convex quadrilateral is $\frac{3C+A}{\pi}$. On the other hand, as already mentioned, another way the points can fail to form the vertices of a convex quadrilateral is if one of them lies inside the triangle formed by the other three. There are four such configurations, depending on which point lies inside the triangle. This gives that the probability that the four points fail to form a quadrilateral is $\frac{4A}{\pi}$. By comparing expectations, we have $3C+A=4A$, so $A=C$ and $4A+3B = \pi$. 
Now consider the expected area of the part of the circle cut off by a chord through two random points $P,Q$ on the side of the chord not containing a third random point, say $R$. This expected area is $B+2C$. We will show that $B+2C = \frac{35}{72\pi} + \frac{\pi}{3}$. Then, using $A=C$ and $4A+3B=\pi$, we get

$4A = \pi - 3(\frac{35}{72\pi} + \frac{\pi}{3} - 2A)$, so $4A =  - \frac{35}{24\pi} + 6A$, i.e. $A = \frac{35}{48\pi}$, as desired.

Now, to compute $B+2C$, let $h$ be the distance from the center $O$ of the circle to the chord $PQ$, and let $h+dh$ be the distance from the center $O$ of the circle to another possible chord $PQ$ where $Q$ can be chosen. The region between the two lines passing through P, with distances $h$ and $h+dh$ from $O$, is the infinitesimal region in which $Q$ can be chosen.

\begin{figure}[h]
  \begin{center}
  \includegraphics[width=2in]{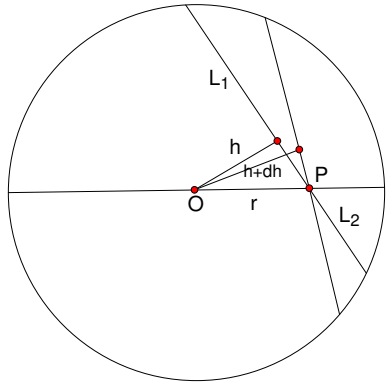}
  \caption{Infinitesimal region}\label{Infinitesimal region}
  \end{center}
\end{figure}

Put $r=OP$. The density of $r$ is $2r dr$ on $[0,1]$. For fixed $r$, $h$ varies over $[0,r]$.
\begin{lem}
The angle between the chords making angles of $\theta = \arcsin(\frac{h}{r})$ and $\alpha = \arcsin(\frac{h+dh}{r})$ with the line containing the point $P$ is $d\theta = \frac{dh}{\sqrt{r^2-h^2}}$.
\end{lem}
\begin{proof}
Let $\theta = \arcsin(\frac{h}{r})$, so
\begin{align*}
\frac{d\theta}{dh} & = \frac{d}{dh} \arcsin\left(\frac{h}{r}\right)\\
& = \frac{1}{\sqrt{1-\left(\frac{h^2}{r^2}\right)}}\frac{1}{r}
\end{align*}
\noindent Hence $d\theta = \frac{dh}{\sqrt{r^2-h^2}}$.
\end{proof}

The infinitesimal region is made up of two sectors, and we can calculate the area of these sectors by calculating the area of the triangles they form. The area of one of the triangles is $\frac{1}{2}L_1 L_1 \sin(d\theta)$, and the other is $\frac{1}{2}L_2
L_2 \sin(d\theta)$, where $L_1, L_2$ are portions of the chord divided by the point $P$. Since the region is so small, $\sin(d\theta)= d\theta$, and so we see that the area of the infinitesimal region which $Q$ can be chosen that we are looking for is $\frac{1}{2}(L_1^2+L_2^2)d\theta$.

\pagebreak

\begin{figure}[h]
  \begin{center}
  \includegraphics[width=2in]{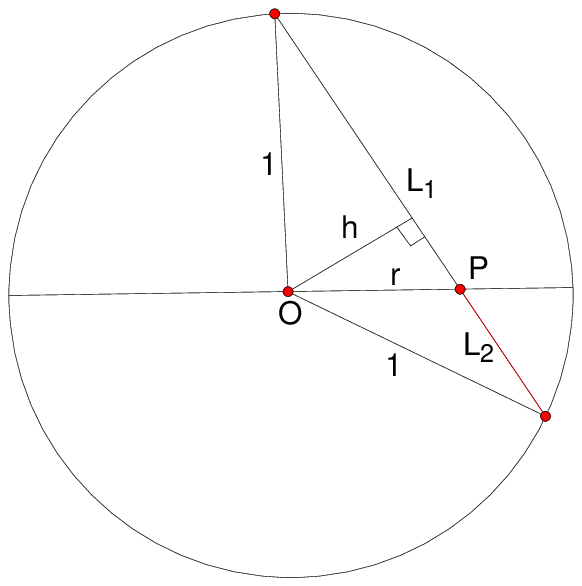}
  \caption{chord PQ}\label{chord PQ}
  \end{center}
\end{figure}

To express the lengths $L_1$ and $L_2$ in terms of $h$ and $r$, we need to look at the three different right triangles in Figure 1.8. Using the Pythagorean theorem, we can find the length of the base of these triangles, and $L_1$ and $L_2$ is made up of some combinations of the lengths of these bases.
Now,
\begin{align*}
L_1 & = \sqrt{1-h^2} + \sqrt{r^2-h^2}\\
L_1^2 & =1-h^2+2 \sqrt{1-h^2} \sqrt{r^2-h^2}+r^2-h^2
\end{align*}
And,
\begin{align*}
L_2 & = \sqrt{1-h^2}-\sqrt{r^2-h^2}\\
L_2^2 & =1-h^2-2 \sqrt{1-h^2} \sqrt{r^2-h^2}+r^2-h^2.
\end{align*}
Thus, $(L_1^2+L_2^2)d\theta =
2(1+r^2-2h^2)\frac{dh}{\sqrt{r^2-h^2}}=\frac{2(1+r^2-2h^2)}{\sqrt{r^2-h^2}}dh$. 
If we divide by $\pi$, and integrate over $r$, the
distribution of $h$ is given by
$$
\frac{1}{\pi}\int_h^1 \frac{2(1+r^2-2h^2)}{\sqrt{r^2-h^2}} 2r dr.
$$
Let $u=r^2-h^2$, so $du=2rdr$, then
\begin{align*}
\frac{1}{\pi}\int_h^1 \frac{2(1+r^2-2h^2)}{\sqrt{r^2-h^2}} 2r dr & = \frac{2}{\pi} \int_0^{1-h^2} \frac{1-h^2+u}{\sqrt{u}}du\\
& = \frac{2}{\pi} (2(1-h^2)^{3/2} + \frac{2}{3}(1-h^2)^{3/2})\\
& =\frac{16}{3\pi}(1-h^2)^{3/2}.
\end{align*}
Let the given unit circle be cut by a chord through two random points $P,Q$. Let $A(h)$ denote the area of the smaller portion, where $h$ is the distance from the center of the circle to the chord. Then
\begin{align*}
A(h) & = 2\int_h^1 \sqrt{1-x^2} dx\\
& =2[\frac{\pi}{4} -  \frac{h}{2} \sqrt{1-h^2} - \frac{1}{2} \arcsin(h)]\\
& =\frac{\pi}{2} -h \sqrt{1-h^2} - \arcsin(h).
\end{align*}

\begin{figure}[h]
  \begin{center}
  \includegraphics[width=2in]{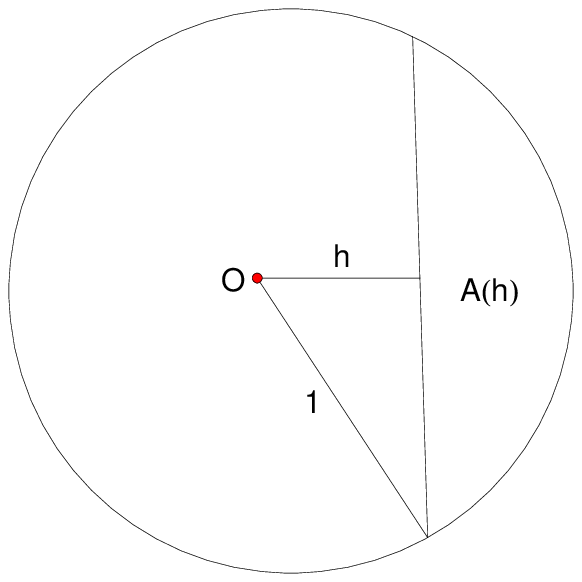}
  \caption{Area A(h)}\label{Area A(h)}
  \end{center}
\end{figure}

Now there are two cases to consider:

Case $a)$: When the third point, $R$, is in the smaller portion, then the area we are
trying to compute is $\pi - A(h)$. This occurs with probability $\frac{A(h)}{\pi}$.

Case $b)$: When the third point, $R$, is in the larger portion, then the area we are
trying to compute is $A(h)$. This occurs with probability $1 - \frac{A(h)}{\pi}$.

Summing these two cases yield $\frac{2}{\pi} A(h)(\pi - A(h))$ and integrating over $h$ yields:
\begin{align*}
B+2C & = \frac{2}{\pi} \int_0^1 A(h)(\pi - A(h)) \frac{16}{3\pi}(1-h^2)^{\frac{3}{2}}dh\\
& =\frac{32}{3\pi^2} \int_0^1 \left(\frac{\pi}{2} -[h \sqrt{1-h^2} + \arcsin(h)]\right)\left(\frac{\pi}{2} + [h \sqrt{1-h^2} + \arcsin(h)]\right) (1-h^2)^{\frac{3}{2}}dh\\
& =\frac{32}{3\pi^2} \int_0^1
\left[\frac{\pi^2}{4}(1-h^2)^{\frac{3}{2}} -h^2(1-h^2)^{\frac{5}{2}}
-2h(1-h^2)^2 \arcsin(h) - \arcsin^2(h)(1-h^2)^{\frac{3}{2}}\right]dh.
\end{align*}
With some calculations, we get that $B+2C = \frac{35}{72\pi} + \frac{\pi}{3}$, (which completes the proof).

\subsection{The third solution}

\noindent The third solution, by Noam Elkies, is summarized below:

Let $P,Q,R$ be three random points inside a unit circle. Let $O$ be the center of the circle, and put $c=$ max$\{OP,OQ,OR\}$. The probability that $OP \leq r$ is the area of the circle with radius $r$ divided by the area of the unit circle. $P(OP \leq r) = \frac{\pi r^2}{\pi} = r^2$. Since the three points are independent, the probability that $c \leq r$ is $(r^2)^3=r^6$. So the density of $c$ is $6c^5 dc$ on $[0,1]$. 
Given $c$, the expectation of $[PQR]$ is equal to $c^2$ times $E$, the expected area of a triangle formed by two random points $A,B$ in a circle and a fixed point $C$ on the boundary.
$$
E([PQR])=c^2E
$$
where $E=E([ABC])$.

We introduce polar coordinates (where $r^2=x^2+y^2, y=r\sin\theta, x=r\cos\theta$) centered at $C$, in which the circle is given by $r=2\sin(\theta)$ for $\theta \in [0,\pi]$.

\begin{figure}[h]
  \begin{center}
  \includegraphics[width=2in]{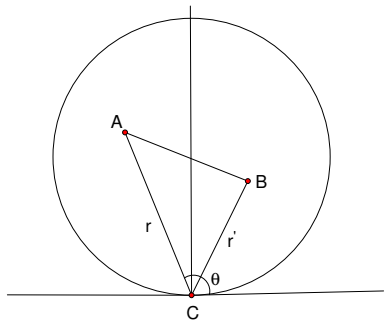}
  \caption{Putnam Third Solution}\label{Putnam Third Solution}
  \end{center}
\end{figure}

The density of a random point in that circle is $\frac{1}{\pi}rdrd\theta$, for $\theta \in [0,\pi]$ and $r \in [0,2\sin(\theta)]$. Let $A(r,\theta)$ and $B(r',\theta')$ be the two random points; then the area of $\triangle ABC$ is $\frac{1}{2}rr'\sin|\theta -
\theta'|$
Since $E=E([ABC])$
\begin{align*}
E & =\int_0^\pi \int_0^\pi \int_0^{2\sin(\theta)} \int_0^{2\sin(\theta')} \frac{1}{2}rr'\sin|\theta - \theta'| \frac{1}{\pi}r \frac{1}{\pi}r' dr'drd\theta'd\theta\\
& = \frac{32}{9\pi^2} \int_0^\pi \int_0^\pi \sin^3(\theta') \sin^3(\theta) \sin|\theta - \theta'| d\theta' d\theta\\
& = \frac{32}{9\pi^2} \left[\int_0^\pi \int_0^\theta
\sin^3(\theta') \sin^3(\theta) \sin(\theta - \theta') d\theta'
d\theta + \int_0^\pi \int_\theta^\pi \sin^3(\theta') \sin^3(\theta)
\sin(\theta' - \theta) d\theta' d\theta\right].
\end{align*}
Let's interchange $\theta$ and $\theta'$ in the second part of the above integral.
$$
E= \frac{32}{9\pi^2} \left[\int_0^\pi \int_0^\theta
\sin^3(\theta') \sin^3(\theta) \sin(\theta - \theta') d\theta'
d\theta + \int_0^\pi \int_{\theta'}^\pi \sin^3(\theta) \sin^3(\theta')
\sin(\theta - \theta') d\theta d\theta'\right].
$$
Therefore $0 \leq \theta' \leq \theta \leq \pi$.
\begin{align*}
E & = \frac{32}{9\pi^2} \left[\int_0^\pi \int_0^\theta \sin^3(\theta') \sin^3(\theta) \sin(\theta - \theta') d\theta' d\theta + \int_0^\theta \int_0^\pi \sin^3(\theta) \sin^3(\theta') \sin(\theta - \theta') d\theta d\theta'\right]\\
& = \frac{64}{9\pi^2} \int_0^\pi \int_0^\theta
\sin^3(\theta') \sin^3(\theta) \sin(\theta - \theta') d\theta'
d\theta.
\end{align*}
With some calculations, we get that $E=\frac{35}{36\pi}$. Therefore
$$
E([PQR]) = \int_0^1 6c^5 c^2 E dc = \int_0^1 6c^7 \frac{35}{36\pi} dc = \frac{35}{48\pi}.
$$

\section{Crofton's Formulas And Proofs}

\noindent We saw in the above arguments that the solution of Sylvester's problem requires the computation of the expected area of a triangle determined by three randomly chosen points in a convex body. Before approaching this problem, we will look at a lower dimensional version, namely, the expected distance between two points chosen at random in a convex body.

\subsection{Crofton's formula}

\noindent Crofton's formula relates the expected number of times a ``random" line intersects a curve in the plane to the length of the curve. 
Given a rectifiable plane curve $\gamma$ (not necessarily closed or simple), and an unoriented line $\sigma$, let $\eta_\gamma (\sigma)$ denote the number of intersection points of $\sigma$ with $\gamma$. This number is finite for almost all $\sigma$ and therefore defines a locally constant function; namely, the value of $\eta_\gamma$ changes when the lines become tangent to the curve $\gamma$. Given an unoriented line, $\sigma$, in the plane, we can find an orthogonal segment from the line to the origin. The line is then determined by the length, $p$, of this segment, and the angle, $\theta$, between the segment and the $x$-axis. Thus $(\theta, p)$ give a coordinate system on the space of unoriented lines, $0 < \theta < 2\pi$, $0 < p < \infty$.

\begin{figure}[h]
  \begin{center}
  \includegraphics[width=2in]{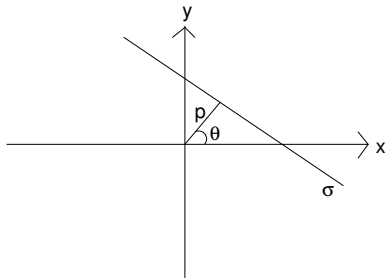}
  \caption{Crofton's formula}\label{Crofton's formula}
  \end{center}
\end{figure}

$m(d\sigma)=d\theta dp$ gives a measure (up to scalars) on this space, and  is isometry invariant. If $(\theta, p)$ are the coordinates of the line $\sigma$, we write $\eta_\gamma(\sigma)=\eta_\gamma(\theta,p)$. We now state Crofton's formula: 
\begin{thm}
Let $\gamma$ be a rectifiable plane curve with length $|\gamma|$. Then $|\gamma| = \frac{1}{2} \int\limits_0^{\infty} \int\limits_0^{2\pi} \eta_\gamma(\theta,p) d\theta dp$
\end{thm}
\begin{proof}
The curve $\gamma$ can be approximated by polygonal lines, and it suffices to prove the theorem for such a line. Suppose that a polygonal line is the concatenation of two, $\gamma_1$ and $\gamma_2$. Both sides of the formula are additive, since $\eta_{\gamma_1 + \gamma_2}(\theta,p) = \eta_{\gamma_1}(\theta,p) + \eta_{\gamma_2}(\theta,p)$, so the formula for $\gamma$ would follow from those for $\gamma_1$ and $\gamma_2$. Hence it suffices to establish the theorem for a segment. This can be done by either a direct computation or, in a more ``lazy" way, as follows. 
Let $\gamma_0$ be the unit segment and let $\int_N \eta_{\gamma_0}(\theta,p) d\theta dp = C$, where $N$ is the space of unoriented lines in the plane. The constant does not depend on the position of the segment because the measure on the space of lines is isometry invariant. Then, by additivity, $\int_N \eta_{\gamma}(\theta,p) d\theta dp = C|\gamma|$ for every segment $\gamma$. And, by the above arguments, $\int_N \eta_{\gamma}(\theta,p) d\theta dp = C|\gamma|$ for every smooth curve $\gamma$. Thus, it remains only to see that $C=2$. Let $\gamma_1$ be the unit circle centered at the origin, so $\eta_{\gamma_1}(\theta,p)=2$ for all $\theta$ and $0 < p < 1$ and zero otherwise. We ignore the case when the line becomes tangent to the curve $\gamma$ because these lines has measure zero. Then,
\begin{align*}
\int \int \eta_\gamma (\theta,p) d\theta dp & = \int_{0}^{1} \int _0^{2\pi}2 d\theta dp\\
& = \int_{0}^{1} 4\pi dp\\
& = 4\pi=2(2\pi).
\end{align*}
Note: $|\gamma| = 2\pi(1) = 2\pi$; therefore, $C=2$. 

As a special case of Crofton's formula, note that if $K \subset \mathbb{R}^2$ is a convex body and $\gamma$ is the boundary of $K$, then
\begin{equation*}
\eta_{\gamma}(\theta,p) =
\begin{cases}
2 & \text{if $\sigma \cap K \neq \emptyset$ and $\sigma$ is not tangent to $\gamma$}\\
1 & \text{if $\sigma$ is tangent to $\gamma$}\\
0 & \text{otherwise}
\end{cases}
\end{equation*}
Thus, Crofton's formula gives $\int\limits_{\sigma\cap K \neq \emptyset} 2 dpd\theta = 2|\gamma|$. In particular, $\int\limits_{\sigma\cap K \neq \emptyset} dpd\theta = |\gamma|$, for convex bodies $K$.

\textbf{Remark:} For each unoriented line, there are two oriented lines, it follows that $\int \int \eta_{\gamma}(\sigma) m(d\sigma) = 4 |\gamma|$ when the integral is taken over oriented lines.
\end{proof}

\subsection{Crofton's Second Theorem}

\noindent We now move on to consider chords of a convex body determined by two randomly chosen points in the body. Given (a pair of) points $P_1(x_1,y_1)$ and $P_2(x_2,y_2)$, we can determine them uniquely by four coordinates $(x_1, y_1, x_2, y_2)$. Then, the uniform density of choosing $P_1$ and $P_2$ at random in a region $K$ of finite area $A$ is
$$
\frac{1}{A^2} dP_1 dP_2 = \frac{1}{A^2} dx_1 dy_1 dx_2 dy_2
$$

\pagebreak

\begin{figure}[h]
  \begin{center}
  \includegraphics[width=2in]{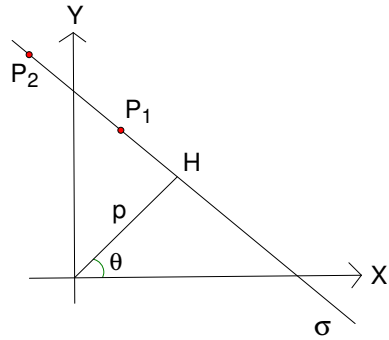}
  \caption{Crofton's Second Theorem}\label{Crofton's Second Theorem}
  \end{center}
\end{figure}

The pair of points can also be determined by the coordinates $(p, \theta)$ of the line $\sigma$ along with the signed distances $t_1, t_2$ from $P_1,P_2$ to the foot $H$ of the perpendicular from the origin $O$ to $\sigma$, see Figure 1.12. To define these signed distances, we need an orientation on $\sigma$, but this can be done by using the right hand rule. 

We wish to write $dP_1 dP_2$ in terms of $p, \theta, t_1, t_2$. From the Figure 1.13, we see that $\cos (\theta) = \frac{x_a}{p}$ and $\cos (\frac{\pi}{2} - \theta) =
\frac{x_b}{t_i}$. Hence $x_a - x_b = x_i = p\cos (\theta) - t_i \sin (\theta)$ for
$i=1,2$. Also, $\sin (\theta) = \frac{y_b}{p}$ and $\sin (\frac{\pi}{2} - \theta) =
\frac{y_a}{t_i}$. Hence $y_a + y_b = y_i = p\sin (\theta) + t_i \cos (\theta)$.

\begin{figure}[h]
  \begin{center}
  \includegraphics[width=2.5in]{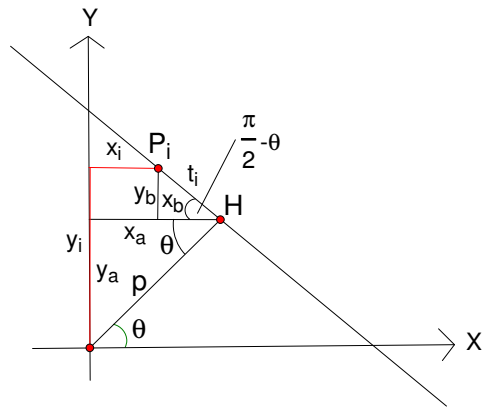}
  \caption{Coordinates of the points}\label{Coordinates of the points}
  \end{center}
\end{figure}

\noindent Thus, taking the derivative we get
\begin{align*}
dx_i & = \cos (\theta) dp - (p \sin (\theta) + t_i \cos (\theta)) d\theta - \sin (\theta) dt_i\\
dy_i & = \sin (\theta) dp + (p \cos (\theta) - t_i \sin (\theta)) d\theta + \cos (\theta) dt_i.
\end{align*}
By exterior multiplication we get,
$$
dx_i \wedge dy_i = p dp \wedge d\theta + dp \wedge dt_i -t_i d\theta \wedge dt_i,
$$
and
$$
dP_1 \wedge dP_2 = dx_1 \wedge dy_1 \wedge dx_2 \wedge dy_2= (t_2 - t_1) dp \wedge d\theta \wedge dt_1 \wedge dt_2.
$$
Since all densities are assumed to be positive, we get
$$
dP_1 dP_2 =  |t_2 - t_1| m(d\sigma) dt_1 dt_2,
$$
where $m(d\sigma) = dp d\theta$ is the usual measure on lines.

Now let $K$ be a closed convex region with area $A$ and boundary curve of length $|\gamma|$. Denote by $l(\sigma)$ the length of the chord determined by the straight line $(p,\theta)$ intersecting $K$.

\begin{figure}[h]
  \begin{center}
  \includegraphics[width=2in]{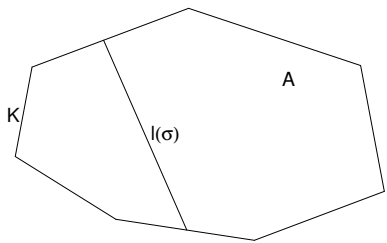}
  \caption{Closed convex region}\label{Closed convex region}
  \end{center}
\end{figure}

Define $I_n = \int\limits_{\sigma \cap K \neq \emptyset} l(\sigma)^n m(d\sigma)$ where $n$ is a positive integer and the integral is taken over all lines $\sigma$ which cut $K$. If $r(P_1,P_2)$ is the  distance between $P_1,P_2$ for any pair of points $P_1,P_2 \in K$, then define $J_n = \int\limits_{P_1,P_2 \in K} r(P_1,P_2)^n dP_1dP_2$. Thus we write $J_n = \int |t_2 - t_1|^{n+1} m(d\sigma)dt_1dt_2$.  
At the end of section 1.2.1, we have noted that it is a consequence of Crofton's formula that  $I_0 = \int m(d\sigma) = |\gamma|$. Also, $I_1 = \int l(\sigma) m(d\sigma) =
\int_0^{\pi} \int l(\sigma) dpd\theta = \pi A$. Crofton's Second Theorem is the following \begin{thm}
$I_3=\int\limits_{\sigma \cap K \neq \emptyset} l(\sigma)^3 m(d\sigma) = 3A^2.$
\end{thm}
\begin{proof}
Fixing $P$, we let $\theta$ determine a chord of length say $l(\sigma)
= b-a$. Fix $t_1$ and integrate $t_2$ over $[a,b]$ and we have, for $n \geq -1$,
\begin{align*}
J_n & = \int m(d\sigma) dt_1 \left[\int_{t_1}^b (t_2 - t_1)^{n+1} dt_2 + \int_a^{t_1} (t_1 - t_2)^{n+1} dt_2\right]\\
& = \frac{1}{n+2} \int m(d\sigma) \int_a^b [(b - t_1)^{n+2} + (t_1 - a)^{n+2}] dt_1\\
& = \frac{1}{n+2} \frac{1}{n+3} \int 2(b - a)^{n+3}m(d\sigma)\\
& = \frac{2}{(n+2)(n+3)} \int l(\sigma)^{n+3}m(d\sigma).
\end{align*}
That is, $J_n = \frac{2}{(n+2)(n+3)} \int l(\sigma)^{n+3}m(d\sigma)$ for $n \geq -1$. Thus, $I_n = \frac{n(n-1)}{2} J_{n-3}$ for $n \geq 2$. Then we have the following:
\begin{align*}
I_2 & = \int l(\sigma)^2 m(d\sigma) = J_{-1} = \int \frac{dP_1dP_2}{r(P_1,P_2)}\\
I_3 & = \int l(\sigma)^3 m(d\sigma) = 3J_0 = 3\int dP_1dP_2 = 3\int dP_1 \int
dP_2 = 3A^2.
\end{align*}
Hence $\int l(\sigma)^3 m(d\sigma) = 3A^2$ for any convex region.
\end{proof}

\section{Distributions And Invariant Measure}

\noindent A measure on a set $S$ is a function that assigns numbers to certain subsets of $S$. Probability is a special kind of measure. The concept has developed in connection with a desire to carry out integration over arbitrary sets. The Lebesgue measure is the standard way of assigning a length, area, or volume to subsets of Euclidean space. 
Lebesgue measure is invariant under rigid motions. Similarly, the space, $SA_{n-1}(\mathbb{R}^n)$, of hyperplanes in $\mathbb{R}^n$ possesses a measure, $m$, which is invariant under rigid motions. To describe this measure, we first need to discuss coordinates for $SA_{n-1}(\mathbb{R}^n)$. If $\sigma \in SA_{n-1}(\mathbb{R}^n)$ is a hyperplane in $\mathbb{R}^n$ which does not pass through the origin, then there is a unique vector, $\vec{v}$, from the origin to $\sigma$ which is orthogonal to $\sigma$. Let $p=|\vec{v}|$; Then we can parametrize points in $SA_{n-1}(\mathbb{R}^n)$ by $(p, \omega)$ where $\omega \in S^{n-1}$. Then the invariant measure on $SA_{n-1}(\mathbb{R}^n)$ is given by $m(d\sigma)=(n-1)! dp$ $m(d\omega)$ where $m(d\omega)$ is the rotationally invariant measure on $S^{n-1}$. We identify this explicitly in coordinates on $S^{n-1}$ for $n=2,3$ and $4$. 

If $n=2$, $m(d\omega)= d\theta$ in terms of the standard angular coordinate, $\theta$ on $S^1$. Hence, $m(d\sigma)=dp d\theta$. Note, this is the measure we used in the previous section. 

For $n=3$, we use standard spherical coordinates, i.e. $x=\cos(\theta)\sin(\phi)$, $y=\sin(\theta)\sin(\phi)$ and $z=\cos(\phi)$ then $m(d\omega)= \sin(\phi)d\theta d\phi$ is the invariant measure on $S^2$. So $m(d\sigma)=2\sin(\phi) dp d\theta d\phi$ is the invariant measure on $SA_2(\mathbb{R}^3)$. 

For $n=4$, we use four-dimensional spherical coordinates, i.e. $x_1=\cos(\theta)\sin(\phi)\sin(\psi)$, $x_2=\sin(\theta)\sin(\phi)\sin(\psi)$, $x_3=\cos(\phi)\sin(\psi)$ and $x_4=\cos(\psi)$ then $m(d\omega)= \sin(\phi)\sin^2(\psi)d\theta d\phi d\psi$ is the invariant measure on $S^3$. So $m(d\sigma)=6\sin(\phi)\sin^2(\psi) dp d\theta d\phi d\psi$ is the invariant measure on $SA_3(\mathbb{R}^4)$.

\chapter{Random Hyperplanes And Sylvester's Problem}

\noindent From section $1.1$, the reader can see that Sylvester's four
point problem can be solved using various methods that deal with
difficult multiple integrals. In this chapter, we will solve
Sylvester's four point problem by looking at random hyperplanes slicing a
convex body. If we look at convex bodies in two dimensions and the linear
secants that intersect them, the results will be Crofton's Second
Theorem. Coordination of convex bodies in higher dimensions will result in a
solution for Sylvester's four point problem in two dimensions and
higher.

\section{The Relationship Among Random Secants, Crofton's Formula And Sylvester's Problem}

\textbf{Linear secants in two dimensions}

Given two points  $P_1$ and $P_2$ in the convex body $K \subset \mathbb{R}^2$. Writing $P_1=(x_1,y_1)$ and $P_2=(x_2,y_2)$, let $X=(x_1,y_1,x_2,y_2) \in \mathbb{R}^4$. Then there exists a unique line $\sigma(X)$ containing these two points. So $dX =dP_1 dP_2=dx_1 dy_1 dx_2 dy_2$.

Recall from section 1.2.2 that
$$
dP_1 dP_2 = |t_2 - t_1|m(d\sigma) dt_1 dt_2
$$
where $m(d\sigma) = dp d\theta$ is the usual measure on lines and $t_1, t_2$ are the signed distances of $P_1, P_2$ to the foot of the perpendicular from the origin to the line $\sigma$.

Now let $\lambda$ denote the probability measure corresponding to the distribution of $\sigma$ given by choosing $P_1$ and $P_2$ uniformly. To find $\lambda$, we have to restrict $dX$ to the subset $K^2=K$ x $K$ and integrate out the $dt_1, dt_2$ over $K\cap\sigma$. So $\lambda(d\sigma) = \frac{1}{A^2}\int\limits_{\sigma} \int\limits_{\sigma} dX$. Let $l(\sigma)$ be the length of $K\cap\sigma$, and let $v(K\cap\sigma)$ denote the normalized expected distance $\left(\frac{E\{|t_2 - t_1|\}}{l(\sigma)}\right)$ between two points on $K\cap\sigma$. Let $A$ be the area of $K$. 
\begin{lem}
Let $2$ points $P_1, P_2$ be chosen at random within the convex body $K \subset \mathbb{R}^2$, and let $\sigma=\sigma(X)$ be the line containing them. Then the distribution $\lambda$ of $\sigma$ is given by
$$
\lambda(d\sigma) = v(K\cap\sigma)l(\sigma)^3A^{-2}m(d\sigma)
$$
where $m(d\sigma)$ is the usual measure on lines.
\end{lem}
\begin{proof}
From
\begin{align*}
v(K\cap\sigma) & =\dfrac{E\{|t_2 - t_1|\}}{l(\sigma)}\\
& =l(\sigma)^{-1}l(\sigma)^{-2}\int\limits_{\sigma}\int\limits_{\sigma}|t_2 - t_1|dt_1 dt_2,
\end{align*}
we infer $\int\limits_{\sigma}\int\limits_{\sigma}|t_2 - t_1|dt_1 dt_2 = v(K\cap\sigma)l(\sigma)^{3}$. Therefore
$$
\lambda(d\sigma) = v(K\cap\sigma)l(\sigma)^3A^{-2}m(d\sigma).
$$
\end{proof}

Let $u$ and $v$ be two points chosen at random from $[0,1]$. The expected value of
the distance between them is
$$
\int_0^1 \int_0^1 |u - v| du dv = 2 \int_0^1 \int_0^{u} (u - v) du dv = \frac{1}{3}.
$$
Thus $v(K\cap\sigma)= \dfrac{1}{3}$ for all $\sigma$.

Let $\sigma(X)$ be the line joining two points chosen at random inside a convex body $K \subset \mathbb{R}^2$ of area $A$. Now, we have that $\lambda(d\sigma) = \frac{1}{3}{l(\sigma)}^{3}A^{-2}m(d\sigma)$. Since $\lambda$ is a probability measure, by definition $\int \lambda (d\sigma) = 1$, and so $\int {l(\sigma)}^{3} m(d\sigma)= 3A^2$.

A result which in two dimensions is Crofton's Second Theorem, as seen in section 1.2.2.
$$
\int\limits_{\sigma\cap K \neq \emptyset}{l(\sigma)}^3 m(d\sigma)=3A^2
$$  

Another, viewpoint to take is this: if $v(K\cap \sigma)$ is not known, but independent of $\sigma$ then by computing $\int {l(\sigma)}^3 m(d\sigma)$ we can find $v(K\cap \sigma)$. This is the view taken to solve Sylvester's problem for the unit ball.

\section{The Solution Of the Sylvester's Problem For A Unit Ball In Two And Three Dimensions}

\noindent In this section, we will apply the results of Kingman's paper in three and four dimensions to solve Sylvester's four point problem in two and three dimensions, respectively.

\subsection{Three Dimensions}

\noindent Now we will solve the Sylvester four point problem in two dimensions by looking at convex bodies in three dimensions and their slices by random hyperplanes.

\noindent \textbf{Planes in three-dimensional space}

Let $K$ be a convex body in $\mathbb{R}^3$, and let $a,b,c$ be three points in $K$. The coordinates of these three points  are $a=(a_1,a_2, a_3)$, $b=(b_1,b_2, b_3)$ and $c=(c_1,c_2, c_3)$. Let $a$, $b$ and $c$ be independently and uniformly distributed over $K$. Now let $X=(a_1,a_2,a_3,b_1,b_2,b_3,c_1,c_2,c_3)$ in $\mathbb{R}^9$. There is a unique plane $\sigma(X)$ in $\mathbb{R}^3$ containing the points $a$, $b$ and $c$, and let $X$ be uniformly distributed over the corresponding subset $K^3$ where $K^3 = \{X \in \mathbb{R}^9 ;a,b,c \in K\}$. Let $\omega$ be the plane parallel to $\sigma(X)$ through the origin; $\omega$ makes an angle $\theta$ with the positive $x_1$-axis, and the unit normal to $\omega$ makes an angle $\phi$ with the positive $x_3$-axis.

\begin{figure}[h]
  \begin{center}
  \includegraphics[width=2in]{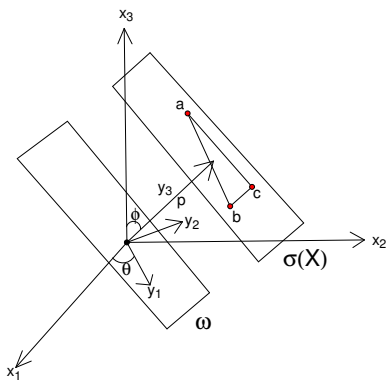}
  \caption{Three-dimensions}\label{Three-dimensions}
  \end{center}
\end{figure}

For any $\omega$, we may choose an orthogonal coordinate system for $\mathbb{R}^3$ in which\\ $\omega = \{y \in \mathbb{R}^3 ; y_3=0\}$. A real square matrix $P$ is orthogonal if $PP^T = P^TP = I$. It is a proper orthogonal matrix if det$(P)=1$. We can use basic trigonometry to find $v_3$ in terms of our new coordinates $\theta, \phi$, where $P=(v_1, v_2, v_3)$ and the columns should be an orthonormal basis. The proper orthogonal matrix connecting the two coordinate systems, so that $x=Py$, is
$$
P=P(\omega)=
\begin{pmatrix}
\cos(\theta)\cos(\phi) & -\sin(\theta) & \cos(\theta)\sin(\phi)\\
\sin(\theta)\cos(\phi) & \cos(\theta) & \sin(\theta)\sin(\phi)\\
-\sin(\phi) & 0 & \cos(\phi)
\end{pmatrix}
$$
Now the plane $\sigma$ is parallel to the unique $\omega$, and $\sigma=\{y \in \mathbb{R}^3 ; y_3= p\}$ where $p \in \mathbb{R}$. We take $p$, $\phi$, and $\theta$ as parametrising the plane $\sigma$. The coordinates $\alpha_1$, $\alpha_2$, $\beta_1$, $\beta_2$, $\gamma_1$, $\gamma_2$ describe the relative positions of the points $a$, $b$, $c$ in the plane $\sigma=\sigma(X)$ with respect to the $Y$ coordinate system. Therefore we have the vector $Y=(\alpha_1,\alpha_2,\beta_1,\beta_2,\gamma_1,\gamma_2, p,\theta, \phi)$, so the mapping $Y \rightarrow X$  is one-to-one. If $\sigma=\sigma(X)$, the point $a$, $b$, $c$ lies in $\sigma$, so that $a=P\alpha$, $b=P\beta$ and $c=P\gamma$, where $\alpha=(\alpha_1,\alpha_2,p)$, $\beta=(\beta_1,\beta_2,p)$, and $\gamma=(\gamma_1,\gamma_2,p)$. Then we get
\begin{align*}
a_1 & =\alpha_1\cos(\theta)\cos(\phi)-\alpha_2\sin(\theta)+p\cos(\theta)\sin(\phi)\\
a_2 & = \alpha_1\sin(\theta)\cos(\phi)+\alpha_2\cos(\theta)+p\sin(\theta)\sin(\phi)\\
a_3 & = -\alpha_1\sin(\phi)+p\cos(\phi)\\
b_1 & = \beta_1\cos(\theta)\cos(\phi)-\beta_2\sin(\theta)+p\cos(\theta)\sin(\phi)\\
b_2 & = \beta_1\sin(\theta)\cos(\phi)+\beta_2\cos(\theta)+p\sin(\theta)\sin(\phi)\\
b_3 & = -\beta_1\sin(\phi)+p\cos(\phi)\\
c_1 & = \gamma_1\cos(\theta)\cos(\phi)-\gamma_2\sin(\theta)+p\cos(\theta)\sin(\phi)\\
c_2 & = \gamma_1\sin(\theta)\cos(\phi)+\gamma_2\cos(\theta)+p\sin(\theta)\sin(\phi)\\
c_3 & = -\gamma_1\sin(\phi)+p\cos(\phi)
\end{align*}
Let $[\alpha, \beta, \gamma]$ be the area of the triangle with vertices $\alpha, \beta, \gamma$ in $\mathbb{R}^2$.
\begin{thm}
Lebesgue measure in X-space is given by
$$
dX = da db dc = [\alpha, \beta, \gamma] d\alpha_{1} d\alpha_{2} d\beta_{1} d\beta_{2} d\gamma_{1} d\gamma_{2} m(d\sigma)
$$
where $m(d\sigma)$ is the invariant measure on $SA_2(\mathbb{R}^3)$.
\end{thm}
\begin{proof} We must compute the Jacobian of the transformation given above to find the distribution of Y. The evaluation of the Jacobian matrix is given by taking the partial derivatives of the coordinates of $X$ with respect to those of $Y$. 
So 
$$
J =
\begin{vmatrix}
\frac{\partial a_1}{\partial \alpha_1} & \frac{\partial a_1}{\partial \alpha_2} & \frac{\partial a_1}{\partial \beta_1} & \frac{\partial a_1}{\partial \beta_2} & \frac{\partial a_1}{\partial \gamma_1} & \frac{\partial a_1}{\partial \gamma_2} & \frac{\partial a_1}{\partial p} & \frac{\partial a_1}{\partial \theta} & \frac{\partial a_1}{\partial \phi}\\
\frac{\partial a_2}{\partial \alpha_1} & \frac{\partial a_2}{\partial \alpha_2} & \frac{\partial a_2}{\partial \beta_1} & \frac{\partial a_2}{\partial \beta_2} & \frac{\partial a_2}{\partial \gamma_1} & \frac{\partial a_2}{\partial \gamma_2} & \frac{\partial a_2}{\partial p} & \frac{\partial a_2}{\partial \theta} & \frac{\partial a_2}{\partial \phi}\\
\frac{\partial a_3}{\partial \alpha_1} & \frac{\partial a_3}{\partial \alpha_2} & \frac{\partial a_3}{\partial \beta_1} & \frac{\partial a_3}{\partial \beta_2} & \frac{\partial a_3}{\partial \gamma_1} & \frac{\partial a_3}{\partial \gamma_2} & \frac{\partial a_3}{\partial p} & \frac{\partial a_3}{\partial \theta} & \frac{\partial a_3}{\partial \phi}\\
\frac{\partial b_1}{\partial \alpha_1} & \frac{\partial b_1}{\partial \alpha_2} & \frac{\partial b_1}{\partial \beta_1} & \frac{\partial b_1}{\partial \beta_2} & \frac{\partial b_1}{\partial \gamma_1} & \frac{\partial b_1}{\partial \gamma_2} & \frac{\partial b_1}{\partial p} & \frac{\partial b_1}{\partial \theta} & \frac{\partial b_1}{\partial \phi}\\
\frac{\partial b_2}{\partial \alpha_1} & \frac{\partial b_2}{\partial \alpha_2} & \frac{\partial b_2}{\partial \beta_1} & \frac{\partial b_2}{\partial \beta_2} & \frac{\partial b_2}{\partial \gamma_1} & \frac{\partial b_2}{\partial \gamma_2} & \frac{\partial b_2}{\partial p} & \frac{\partial b_2}{\partial \theta} & \frac{\partial b_2}{\partial \phi}\\
\frac{\partial b_3}{\partial \alpha_1} & \frac{\partial b_3}{\partial \alpha_2} & \frac{\partial b_3}{\partial \beta_1} & \frac{\partial b_3}{\partial \beta_2} & \frac{\partial b_3}{\partial \gamma_1} & \frac{\partial b_3}{\partial \gamma_2} & \frac{\partial b_3}{\partial p} & \frac{\partial b_3}{\partial \theta} & \frac{\partial b_3}{\partial \phi}\\
\frac{\partial c_1}{\partial \alpha_1} & \frac{\partial c_1}{\partial \alpha_2} & \frac{\partial c_1}{\partial \beta_1} & \frac{\partial c_1}{\partial \beta_2} & \frac{\partial c_1}{\partial \gamma_1} & \frac{\partial c_1}{\partial \gamma_2} & \frac{\partial c_1}{\partial p} & \frac{\partial c_1}{\partial \theta} & \frac{\partial c_1}{\partial \phi}\\
\frac{\partial c_2}{\partial \alpha_1} & \frac{\partial c_2}{\partial \alpha_2} & \frac{\partial c_2}{\partial \beta_1} & \frac{\partial c_2}{\partial \beta_2} & \frac{\partial c_2}{\partial \gamma_1} & \frac{\partial c_2}{\partial \gamma_2} & \frac{\partial c_2}{\partial p} & \frac{\partial c_2}{\partial \theta} & \frac{\partial c_2}{\partial \phi}\\
\frac{\partial c_3}{\partial \alpha_1} & \frac{\partial c_3}{\partial \alpha_2} & \frac{\partial c_3}{\partial \beta_1} & \frac{\partial c_3}{\partial \beta_2} & \frac{\partial c_3}{\partial \gamma_1} & \frac{\partial c_3}{\partial \gamma_2} & \frac{\partial c_3}{\partial p} & \frac{\partial c_3}{\partial \theta} & \frac{\partial c_3}{\partial \phi}
\end{vmatrix}
.
$$
This gives $J =$
$
\begin{vmatrix}
A & 0 & 0 & B(\alpha)\\
0 & A & 0 & B(\beta)\\
0 & 0 & A & B(\gamma)
\end{vmatrix}
$, where $A=$
$
\begin{pmatrix}
\cos(\theta)\cos(\phi) & -\sin(\theta)\\
\sin(\theta)\cos(\phi) & \cos(\theta)\\
-\sin(\phi) & 0
\end{pmatrix}
$\\

\noindent and $B(x)$ is the matrix

$
\begin{small}
\begin{pmatrix}
\cos(\theta)\sin(\phi) & -x_1\sin(\theta)\cos(\phi)- x_2\cos(\theta)-p\sin(\theta)\sin(\phi) & -x_1\cos(\theta)\sin(\phi)+p\cos(\theta)\cos(\phi)\\
\sin(\theta)\sin(\phi) & x_1\cos(\theta)\cos(\phi)- x_2\sin(\theta)+p\cos(\theta)\sin(\phi) & -x_1\sin(\theta)\sin(\phi)+p\sin(\theta)\cos(\phi)\\
\cos(\phi) & 0 & -x_1\cos(\phi)-p\sin(\phi)
\end{pmatrix}
\end{small}
$.\\

Note also that the $0$ in the above matrix $J$ is a $3$ x $2$ zero matrix. To find J, we will multiply its matrix by three copies of the transpose of P, and take its determinant. That is, let
$$
R =
\begin{pmatrix}
P^T & 0 & 0 \\
0 & P^T & 0\\
0 & 0 & P^T
\end{pmatrix}
\begin{pmatrix}
A & 0 & 0 & B(\alpha)\\
0 & A & 0 & B(\beta)\\
0 & 0 & A & B(\gamma)
\end{pmatrix},
$$ 
and note $J=|R|$. 
Note the $0$ in the first matrix above is a $3$ x $3$ zero matrix, and the $0$ in the second matrix above is a $3$ x $2$ zero matrix. So the product is given by
$$
R =
\begin{pmatrix}
1 & 0 & 0 & 0 & 0 & 0 & 0 & -\alpha_2\cos(\phi) & p\\
0 & 1 & 0 & 0 & 0 & 0 & 0 & \alpha_1\cos(\phi)+p\sin(\phi) & 0\\
0 & 0 & 0 & 0 & 0 & 0 & 1 & -\alpha_2\sin(\phi) & -\alpha_1\\
0 & 0 & 1 & 0 & 0 & 0 & 0 & -\beta_2\cos(\phi) & p\\
0 & 0 & 0 & 1 & 0 & 0 & 0 & \beta_1\cos(\phi)+p\sin(\phi) & 0\\
0 & 0 & 0 & 0 & 0 & 0 & 1 & -\beta_2\sin(\phi) & -\beta_1\\
0 & 0 & 0 & 0 & 1 & 0 & 0 & -\gamma_2\cos(\phi) & p\\
0 & 0 & 0 & 0 & 0 & 1 & 0 & \gamma_1\cos(\phi)+p\sin(\phi) & 0\\
0 & 0 & 0 & 0 & 0 & 0 & 1 & -\gamma_2\sin(\phi) & -\gamma_1
\end{pmatrix}
.
$$
So $J=|R|=(1)(1)(-1)(-1)(1)(1)
\begin{vmatrix}
1 & -\alpha_2\sin(\phi) & -\alpha_1\\
1 & -\beta_2\sin(\phi) & -\beta_1\\
1 & -\gamma_2\sin(\phi) & -\gamma_1
\end{vmatrix}
$, giving
$$
J=
\begin{vmatrix}
1 & -\alpha_2\sin(\phi) & -\alpha_1\\
0 & (\alpha_2-\beta_2)\sin(\phi) & \alpha_1-\beta_1\\
0 & (\alpha_2-\gamma_2)\sin(\phi) & \alpha_1-\gamma_1
\end{vmatrix}
,
$$
so 
$$
J=
\begin{vmatrix}
(\alpha_2-\beta_2)\sin(\phi) & \alpha_1-\beta_1\\
(\alpha_2-\gamma_2)\sin(\phi) & \alpha_1-\gamma_1
\end{vmatrix}
=\sin(\phi)
\begin{vmatrix}
\alpha_2-\beta_2 & \alpha_1-\beta_1\\
\alpha_2-\gamma_2 & \alpha_1-\gamma_1
\end{vmatrix}
.
$$
We have $dX = J dY$, so
\begin{align*}
dX=dadbdc & =\sin(\phi)\begin{vmatrix}
\alpha_2-\beta_2 & \alpha_1-\beta_1\\
\alpha_2-\gamma_2 & \alpha_1-\gamma_1
\end{vmatrix}
d\alpha_1 d\alpha_2 d\beta_1 d\beta_2 d\gamma_1 d\gamma_2 dp d\theta d\phi\\
& = \frac{1}{2}\begin{vmatrix}
\alpha_2-\beta_2 & \alpha_1-\beta_1\\
\alpha_2-\gamma_2 & \alpha_1-\gamma_1
\end{vmatrix}
d\alpha_{1} d\alpha_2 d\beta_1 d\beta_2 d\gamma_1 d\gamma_2 m(d\sigma)
\end{align*}
where $m(d\sigma) = 2dp d\theta d\phi \sin(\phi)$ is the invariant measure on $SA_2(\mathbb{R}^3)$.

Now, the area of a triangle with vertices $\alpha, \beta, \gamma$ is
$$
[\alpha, \beta, \gamma] = \frac{1}{2}\left\bracevert\begin{vmatrix}
\beta_1-\alpha_1 & \beta_2-\alpha_2\\
\gamma_1-\alpha_1 & \gamma_2-\alpha_2
\end{vmatrix}\right\bracevert
$$
Note the double line means the absolute value of the determinant. Thus we get that $dX = [\alpha, \beta, \gamma] d\alpha_{1} d\alpha_{2} d\beta_{1} d\beta_{2} d\gamma_{1} d\gamma_{2} m(d\sigma)$.
\end{proof}

To find $\lambda$, we have to restrict $dX$ to the subset $K^3$ and integrate out the $d\alpha_1, d\alpha_2, d\beta_1, d\beta_2, d\gamma_1, d\gamma_2$ over $K\cap\sigma$. Let $A_\sigma$ denote the area of $K\cap\sigma$, and $v(K\cap\sigma)=A_{\sigma}^{-1}E[\alpha,\beta,\gamma]$. Also, let $V$ denote the volume of $K$. Then
\begin{lem}
Let $3$ points $a,b,c$ be chosen at random within the convex body $K \subset \mathbb{R}^3$, and let $\sigma$ be the plane containing them. Then the distribution $\lambda$ of $\sigma$ is given by
$$
\lambda(d\sigma) = v(K\cap\sigma){A_\sigma}^4V^{-3}m(d\sigma)
$$
where $m(d\sigma)$ is the invariant measure on $SA_2(\mathbb{R}^3)$.
\end{lem}
\begin{proof} Since $v(K\cap\sigma)=\dfrac{E[\alpha,\beta,\gamma]}{A_\sigma}$, we have 
$$
v(K\cap\sigma)={A_\sigma}^{-1}{A_\sigma}^{-3}\int\limits_{K\bigcap\sigma}\cdots\int\limits_{K\cap\sigma}[\alpha,\beta,\gamma]d\alpha_1d\alpha_2d\beta_1 d\beta_2d\gamma_1d\gamma_2.
$$
Thus,
\begin{align*}
\int\limits_{K\cap\sigma}\cdots\int\limits_{K\cap\sigma}\frac{1}{2}\left\bracevert\begin{vmatrix}
\beta_1-\alpha_1 & \beta_2-\alpha_2\\
\gamma_1-\alpha_1 & \gamma_2-\alpha_2
\end{vmatrix}\right\bracevert d\alpha_1d\alpha_2d\beta_1 d\beta_2d\gamma_1d\gamma_2 = v(K\cap\sigma){A_\sigma}^{4}
\end{align*}
Using the fact that $\int_{K^3}\lambda(d\sigma)=1$, but $\int_{K^3}dX=V^3$, so we get
$$
\lambda(d\sigma) = v(K\bigcap\sigma){A_\sigma}^4V^{-3}m(d\sigma).
$$
\end{proof}

From section 1.3, we have $m(d\sigma) = 2\sin(\phi)dp d\theta d\phi$, and  $m(d\omega) = \sin(\phi) d\theta d\phi$. So $m(d\sigma) = 2dp$ $m(d\omega)$, and we get that
$$
\lambda(d\sigma) = 2v(K\cap\sigma){A_\sigma}^{4}V^{-3}dp m(d\omega)
$$
Referring back to $\int_{K^3}\lambda(d\sigma)=1$, we get
$$
V^{3}= 2\int_{\Sigma}v(K\cap\sigma){A_\sigma}^{4}dp m(d\omega).
$$
where $\Sigma$ is the subset of the space of planes in $\mathbb{R}^3$ which intersect $K$. 

Now we consider the case where $K=B^3$, the unit ball in $\mathbb{R}^3$. It has volume $\beta_3=\frac{4}{3}\pi$ and boundary $\partial K = S^{2}$ has surface area $3\beta_3=4\pi$. Note $K\cap\sigma \approx B^{2}$ with radius $(1-p^2)^{\frac{1}{2}}$, and $B^2$ has area $(1-p^2)\pi$. Also, $v(K\cap\sigma)=v(B^2)$ is independent of $\sigma$, so it can be taken out of the integral.
So we get:  $\beta_3^3 = 2\int_{\Sigma}v(B^{2})((1-p^2)\pi)^{4}dp$ $m(d\omega)$
\begin{align*}
\left(\frac{4}{3}\pi\right)^3 & = 2\int_{\Sigma}v(B^{2})(1-p^2)^4\pi ^{4}dp m(d\omega)\\
& = 2v(B^{2})4\pi(\pi^{4})\int_0^1(1-p^2)^{4}dp.
\end{align*}
Also,
\begin{align*}
\int_0^1(1-p^2)^{4}dp & = \int_0^1 (1-4p^2+6p^4-4p^6+p^8)dp\\
& = \left[p-\frac{4}{3}p^3+\frac{6}{5}p^5-\frac{4}{7}p^7+\frac{1}{9}p^9\right]_0^1\\
& = 1-\frac{4}{3}+\frac{6}{5}-\frac{4}{7}+\frac{1}{9}\\
& = \dfrac{128}{315}.
\end{align*}

We now solve for $v(B^{2})$.
$$
\left(\frac{4}{3}\pi\right)^3 = 2v(B^{2})4\pi\left(\pi ^{4}\right)\frac{128}{315}.
$$
Therefore $v(B^2)=\dfrac{35}{48\pi^2}$.

\begin{center}
\textbf{THE PROBLEM OF SYLVESTER FOR A UNIT BALL}
\end{center}

\noindent \textbf{THE PROBLEM:} Find the probability that four
points chosen at random in the interior of the unit circle form a
convex quadrilateral, that is, that none of the points is inside the
triangle formed by the other three.

\noindent \textbf{THE SOLUTION IN TWO DIMENSIONS:}

The only way the quadrilateral can fail to be convex is for one of the points to lie inside the triangle formed by the other three points. There are four such configurations, depending on which point lies inside the triangle, and they are mutually exclusive, because two points cannot both lie
inside the triangle at the same time. We may ignore the case where three of the points are collinear, as this occurs with probability zero. The probability
that a point lies inside the triangle is given by $v(B^2)$.

\textbf{Answer to Sylvester Problem for a unit circle =} $1-4\left(\dfrac{35}{48\pi^2}\right)= \mathbf{1-\dfrac{35}{12\pi^2}}$

\subsection{Four Dimensions}

\noindent Now we will solve the Sylvester four point problem in three dimensions by looking at convex bodies in four dimensions and their slices by random hyperplanes.

\noindent \textbf{Hyperplanes cutting a sphere}

Let $K$ be a convex body in $\mathbb{R}^4$, and let $a,b,c,d$ be four points in $K$. The coordinates of these four points are $a=(a_1,a_2, a_3,a_4)$, $b=(b_1,b_2, b_3,b_4)$, $c=(c_1,c_2, c_3,c_4)$ and $d=(d_1,d_2, d_3,d_4)$. Let $a$, $b$, $c$ and $d$ be independently and uniformly distributed over $K$. Now let $X=(a_1,a_2,a_3,a_4,b_1,b_2,b_3, b_4,c_1,c_2,c_3,c_4,d_1,d_2, d_3,d_4)$ in $\mathbb{R}^{16}$. There is a unique $3$-flat (hyperplane) $\sigma(X)$ in $\mathbb{R}^4$ containing the points $a$, $b$, $c$ and $d$, and let $X$ be uniformly distributed over the corresponding subset $K^4$ where $K^4 = \{X \in \mathbb{R}^{16}; a,b,c,d \in K\}$. Let $\omega$ be the hyperplane parallel to $\sigma(X)$ through the origin, the unit normal to $\omega$ makes an angle $\psi$ with the positive $x_4$-axis and $\omega$ makes angles $\theta$, and $\phi$ with the $x_1$, $x_2$, $x_3$ plane.
For any $\omega$, we may choose an orthogonal coordinate system for $\mathbb{R}^4$ in which $\omega = \{y \in \mathbb{R}^4 ; y_4=0\}$. The proper orthogonal matrix connecting the two coordinate systems, so that $x=Py$, is
$$
P=P(\omega)=
\begin{pmatrix}
\cos(\theta)\cos(\phi) & -\sin(\theta) & \cos(\theta)\sin(\phi)\cos(\psi) & \cos(\theta)\sin(\phi)\sin(\psi)\\
\sin(\theta)\cos(\phi) & \cos(\theta) & \sin(\theta)\sin(\phi)\cos(\psi) & \sin(\theta)\sin(\phi)\sin(\psi)\\
-\sin(\phi) & 0 & \cos(\phi)\cos(\psi) & \cos(\phi)\sin(\psi)\\
0 & 0 & -\sin(\psi) & \cos(\psi)
\end{pmatrix}
$$
Now the hyperplane $\sigma$ is parallel to the unique $\omega$, and $\sigma=\{y \in \mathbb{R}^4 ; y_4= p\}$ where $p \in \mathbb{R}$. We take $p$, $\phi$, $\theta$ and $\psi$ as parametrising the hyperplane $\sigma$. The coordinates $\alpha_1$, $\alpha_2$, $\alpha_3$, $\beta_1$, $\beta_2$, $\beta_3$, $\gamma_1$, $\gamma_2$, $\gamma_3$, $\delta_1$, $\delta_2$, $\delta_3$  describe the relative positions of the points $a$, $b$, $c$, $d$ in the hyperplane $\sigma=\sigma(X)$ with respect to the $Y$ coordinate system. Therefore we have the vector $Y=(\alpha_1,\alpha_2,\alpha_3,\beta_1,\beta_2,\beta_3,\gamma_1,\gamma_2,\gamma_3,\delta_1,\delta_2,
\delta_3,p,\theta, \phi, \psi)$, so the mapping $Y \rightarrow X$ is one-to-one. If $\sigma=\sigma(X)$, the point $a$, $b$, $c$, $d$ lies in $\sigma$, so that $a=P\alpha$, $b=P\beta$, $c=P\gamma$ and $d=P\delta$, where $\alpha=(\alpha_1,\alpha_2,\alpha_3,p)$,
$\beta=(\beta_1,\beta_2,\beta_3,p)$, $\gamma=(\gamma_1,\gamma_2,\gamma_3,p)$, and
$\delta=(\delta_1,\delta_2,\delta_3,p)$. Then we get
\begin{align*}
a_1 & = \alpha_1\cos(\theta)\cos(\phi)-\alpha_2\sin(\theta)+\alpha_3\cos(\theta)\sin(\phi)\cos(\psi)+p\cos(\theta)\sin(\phi)\sin(\psi)\\
a_2 & = \alpha_1\sin(\theta)\cos(\phi)+\alpha_2\cos(\theta)+\alpha_3\sin(\theta)\sin(\phi)\cos(\psi)+p\sin(\theta)\sin(\phi)\sin(\psi)\\
a_3 & = -\alpha_1\sin(\phi)+\alpha_3\cos(\phi)\cos(\psi)+p\cos(\phi)\sin(\psi)\\
a_4 & = -\alpha_3\sin(\psi)+p\cos(\psi)\\
b_1 & = \beta_1\cos(\theta)\cos(\phi)-\beta_2\sin(\theta)+\beta_3\cos(\theta)\sin(\phi)\cos(\psi)+p\cos(\theta)\sin(\phi)\sin(\psi)\\
b_2 & = \beta_1\sin(\theta)\cos(\phi)+\beta_2\cos(\theta)+\beta_3\sin(\theta)\sin(\phi)\cos(\psi)+p\sin(\theta)\sin(\phi)\sin(\psi)\\
b_3 & = -\beta_1\sin(\phi)-\beta_3\cos(\phi)\cos(\psi)+p\cos(\phi)\sin(\psi)\\
b_4 & = -\beta_3\sin(\psi)+p\cos(\psi)\\
c_1 & = \gamma_1\cos(\theta)\cos(\phi)-\gamma_2\sin(\theta)+\gamma_3\cos(\theta)\sin(\phi)\cos(\psi)+p\cos(\theta)\sin(\phi)\sin(\psi)\\
c_2 & = \gamma_1\sin(\theta)\cos(\phi)+\gamma_2\cos(\theta)+\gamma_3\sin(\theta)\sin(\phi)\cos(\psi)+p\sin(\theta)\sin(\phi)\sin(\psi)\\
c_3 & = -\gamma_1\sin(\phi)+\gamma_3\cos(\phi)\cos(\psi)+p\cos(\phi)\sin(\psi)\\
c_4 & = -\gamma_3\sin(\psi)+p\cos(\psi)\\
d_1 & = \delta_1\cos(\theta)\cos(\phi)-\delta_2\sin(\theta)+\delta_3\cos(\theta)\sin(\phi)\cos(\psi)+p\cos(\theta)\sin(\phi)\sin(\psi)\\
d_2 & = \delta_1\sin(\theta)\cos(\phi)+\delta_2\cos(\theta)+\delta_3\sin(\theta)\sin(\phi)\cos(\psi)+p\sin(\theta)\sin(\phi)\sin(\psi)\\
d_3 & = -\delta_1\sin(\phi)+\delta_3\cos(\phi)\cos(\psi)+p\cos(\phi)\sin(\psi)\\
d_4 & = -\delta_3\sin(\psi)+p\cos(\psi)
\end{align*}
Let $\Delta(\alpha, \beta, \gamma, \delta)$ be the volume of the tetrahedron with vertices $\alpha, \beta, \gamma, \delta \in \mathbb{R}^3$.
\begin{thm}
Lebesgue measure in X-space is given by
$$
dX = da db dc dd = \Delta(\alpha, \beta, \gamma, \delta) d\alpha_{1} d\alpha_{2} d\alpha_{3} d\beta_{1} d\beta_{2} d\beta_{3} d\gamma_{1} d\gamma_{2} d\gamma_{3} d\delta_1 d\delta_2 d\delta_3 m(d\sigma)
$$
where $m(d\sigma)$ is the invariant measure on $SA_3(\mathbb{R}^4)$.
\end{thm}
\begin{proof} We must compute the Jacobian of the transformation given above to find the distribution of Y. The evaluation of the Jacobian Matrix is given by taking the partial derivatives of the coordinates of $X$ with respect to those of $Y$. So 
$$
J =
\begin{vmatrix}
A & 0 & 0 & 0 & B_1(\alpha) & B_2(\alpha)\\
0 & A & 0 & 0 & B_1(\beta) & B_2(\beta)\\
0 & 0 & A & 0 & B_1(\gamma) & B_2(\gamma)\\
0 & 0 & 0 & A & B_1(\delta) & B_2(\delta)
\end{vmatrix}
$$
where $A=$
$
\begin{pmatrix}
\cos(\theta)\cos(\phi) & -\sin(\theta) & \cos(\theta)\sin(\phi)\cos(\psi)\\
\sin(\theta)\cos(\phi) & \cos(\theta) & \sin(\theta)\sin(\phi)\cos(\psi)\\
-\sin(\phi) & 0 & \cos(\phi)\cos(\psi)\\
0 & 0 & -\sin(\psi)
\end{pmatrix}
$\\
 
and, $B_1(x)=$ is the matrix

$
\begin{small}
\begin{pmatrix}
\cos(\theta)\sin(\phi)\sin(\psi) &  -x_1\sin(\theta)\cos(\phi) - x_2\cos(\theta)-x_3\sin(\theta)\sin(\phi)\cos(\psi)-p\sin(\theta)\sin(\phi)\sin(\psi)\\
\sin(\theta)\sin(\phi)\sin(\psi) & x_1\cos(\theta)\cos(\phi) - x_2\sin(\theta) +x_3\cos(\theta)\sin(\phi)\cos(\psi)+p\cos(\theta)\sin(\phi)\sin(\psi)\\
\cos(\phi)\sin(\psi) & 0\\
\cos(\psi) & 0
\end{pmatrix}
\end{small}
$\\ 

and, $B_2(x)=$ is the matrix

$
\begin{scriptsize}
\begin{pmatrix}
-x_1\cos(\theta)\sin(\phi)+x_3\cos(\theta)\cos(\phi)\cos(\psi)+p\cos(\theta)\cos(\phi)\sin(\psi) &  -x_3\cos(\theta)\sin(\phi)\sin(\psi) +p\cos(\theta)\sin(\phi)\cos(\psi)\\
-x_1\sin(\theta)\sin(\phi)+x_3\sin(\theta)\cos(\phi)\cos(\psi)+p\sin(\theta)\cos(\phi)\sin(\psi) &-x_3\sin(\theta)\sin(\phi)\sin(\psi) +p\sin(\theta)\sin(\phi)\cos(\psi)\\
-x_1\cos(\phi) -x_3\sin(\phi)\cos(\psi) -p\sin(\phi)\sin(\psi) & -x_3\cos(\phi)\sin(\psi)+p\cos(\phi)\cos(\psi)\\
0 & -x_3\cos(\psi)-p\sin(\psi)
\end{pmatrix}
\end{scriptsize}
$.\\

Note also that the $0$ in the above matrix $J$ is a $4$ x $3$ zero matrix. To find J, we will multiply its matrix by three copies of the transpose of P, and take its determinant. That is, let
$$
R =
\begin{pmatrix}
P^T & 0 & 0 & 0\\
0 & P^T & 0 & 0\\
0 & 0 & P^T & 0\\
0 & 0 & 0 & P^T
\end{pmatrix}
\begin{pmatrix}
A & 0 & 0 & 0 & B(\alpha)\\
0 & A & 0 & 0 & B(\beta)\\
0 & 0 & A & 0 & B(\gamma)\\
0 & 0 & 0 & A & B(\delta)
\end{pmatrix}
,
$$
and note $J=|R|$. Note the $0$ in the first matrix above is a $4$ x $4$ zero matrix, and the $0$ in the second matrix above is a $4$ x $3$ zero matrix. So the product gives us a $16$ x $16$ matrix with zeros and ones in the first thirteen columns. Doing the same process as in the three-dimensional case, we get that
$J = \sin(\phi)\sin^2(\psi)\begin{vmatrix}
\beta_2 - \alpha_2 & \beta_1 - \alpha_1 & \alpha_3-\beta_3\\
\gamma_2 - \alpha_2 & \gamma_1 - \alpha_1 & \alpha_3-\gamma_3\\
\delta_2 - \alpha_2 & \delta_1 - \alpha_1 & \alpha_3-\delta_3
\end{vmatrix}$.

\noindent We have $dX =JdY$, so
$$
dX = da db dc dd = \frac{1}{6}\begin{vmatrix}
\beta_2 - \alpha_2 & \beta_1 - \alpha_1 & \alpha_3-\beta_3\\
\gamma_2 - \alpha_2 & \gamma_1 - \alpha_1 & \alpha_3-\gamma_3\\
\delta_2 - \alpha_2 & \delta_1 - \alpha_1 & \alpha_3-\delta_3
\end{vmatrix} d\alpha_1 d\alpha_2 d\alpha_3 d\beta_1 d\beta_2 d\beta_3
d\gamma_1 d\gamma_2 d\gamma_3 d\delta_1 d\delta_2 d\delta_3
m(d\sigma)
$$
where $m(d\sigma) = 6dp d\theta d\phi d\psi \sin(\phi)\sin^2(\psi)$ is the invariant measure on $SA_3(\mathbb{R}^4)$.

In the tetrahedron, we will translate the point $(\alpha_1,\alpha_2,\alpha_3)$ to the origin. So, the volume of the tetrahedron with vertices $\alpha, \beta, \gamma, \delta$ is
$$
\Delta(\alpha, \beta, \gamma, \delta)  = \frac{1}{6}
\left\bracevert\begin{vmatrix}
\beta_1-\alpha_1 & \gamma_1-\alpha_1 & \delta_1-\alpha_1\\
\beta_2-\alpha_2 & \gamma_2-\alpha_2 & \delta_2-\alpha_2\\
\beta_3-\alpha_3 & \gamma_3-\alpha_3 & \delta_3-\alpha_3
\end{vmatrix}\right\bracevert
$$
Thus we get $dX = \Delta(\alpha, \beta, \gamma, \delta) d\alpha_{1} d\alpha_{2} d\alpha_{3} d\beta_{1} d\beta_{2} d\beta_{3} d\gamma_{1} d\gamma_{2} d\gamma_{3} d\delta_1 d\delta_2 d\delta_3 m(d\sigma)$.
\end{proof}

To find $\lambda$, we have to restrict $dX$ to the subset $K^4$ and integrate out the $d\alpha_1,d\alpha_2, d\alpha_3, d\beta_1, d\beta_2, d\beta_3, d\gamma_1, d\gamma_2, d\gamma_3, d\delta_1, d\delta_2, d\delta_3$ over $K\cap\sigma$. Let $V_{\sigma}$ denote the $3-$dimensional volume of $K\cap\sigma$, and $v(K\cap\sigma)=V_{\sigma}^{-1}E\{\Delta(\alpha,\beta,\gamma,\delta)\}$. Also, let $V$ denote the $4$-dimensional volume of $K$. Then
\begin{lem}
Let $4$ points $a,b,c,d$ be chosen at random within the convex body $K \subset \mathbb{R}^4$, and let $\sigma$ be the $3$-flat containing them. Then the distribution $\lambda$ of $\sigma$ is given by
$$
\lambda(d\sigma) = v(K\cap\sigma)V_{\sigma}^5V^{-4}m(d\sigma)
$$
where $m$ is the invariant measure on $SA_3(\mathbb{R}^4)$.
\end{lem}
\begin{proof} Since $v(K\cap\sigma)=\dfrac{E\{\Delta(\alpha,\beta,\gamma,\delta)\}}{V_{\sigma}}$, we have
\begin{align*}
v(K\cap\sigma)& = V_{\sigma}^{-1}V_{\sigma}^{-4}\int\limits_{K\cap\sigma}\cdots\int\limits_{K\cap\sigma}\Delta(\alpha,\beta,\gamma,\delta)d\alpha_1d\alpha_2 d\alpha_3d\beta_1d\beta_2d\beta_3d\gamma_1d\gamma_2d\gamma_3d\delta_1d\delta_2d\delta_3\\
v(K\cap\sigma)V_{\sigma}^{5} & = \int\limits_{K\cap\sigma}\cdots\int\limits_{K\cap\sigma}\frac{1}{6}\left\bracevert\begin{vmatrix}
\beta_2 - \alpha_2 & \beta_1 - \alpha_1 & \alpha_3-\beta_3\\
\gamma_2 - \alpha_2 & \gamma_1 - \alpha_1 & \alpha_3-\gamma_3\\
\delta_2 - \alpha_2 & \delta_1 - \alpha_1 & \alpha_3-\delta_3
\end{vmatrix}\right\bracevert d\alpha_1d\alpha_2 d\alpha_3d\beta_1d\beta_2d\beta_3d\gamma_1d\gamma_2d\gamma_3d\delta_1d\delta_2d\delta_3
\end{align*}
Again, using the fact that $\int_{K^4}\lambda(d\sigma)=1$, but $\int_{K^4} dX =V^4$, so we get
$$
\lambda(d\sigma) = v(K\cap\sigma)V_{\sigma}^5V^{-4}m(d\sigma).
$$
\end{proof}

From section 1.3, we have $m(d\sigma) = 6\sin(\phi)\sin^2(\psi)dp d\theta d\phi
d\psi$, and $m(d\omega) = \sin(\phi)\sin^2(\psi) d\theta d\phi d\psi$. So $m(d\sigma) = 3!dp$ $m(d\omega)$, and we get that
$$
\lambda(d\sigma) = 3!v(K\cap\sigma)V_{\sigma}^{5}V^{-4}dpm(d\omega)
$$
Referring back to $\int_{K^4}\lambda(d\sigma)=1$, we get
$$
V^{4}= 3!\int_{\Sigma}v(K\cap\sigma)V_{\sigma}^{5}dpm(d\omega).
$$
where $\Sigma$ is the subset of the space of 3-flats in $\mathbb{R}^4$ which intersect $K$. 

Now we consider the case where $K=B^4$, the unit ball in $\mathbb{R}^4$. It has volume $\beta_4$ and $\partial K = S^{3}$ has surface area $4\beta_4$. Note $K\cap\sigma \approx B^{3}$ with radius $(1-p^2)^{\frac{1}{2}}$, and $B^{3}$  has volume $(1-p^2)^{\frac{3}{2}}\beta_{3}$. Also, $v(K\cap\sigma)=v(B^3)$ is independent of $\sigma$, so it can be taken out of the integral. 
So we get:  $\beta_4^4 =
3!\int_{\Sigma}v(B^{3})((1-p^2)^{{\frac{3}{2}}}\beta_{3})^{5}dpm(d\omega)$
\begin{align*}
\beta_4^4 & = 3!\int_{\Sigma}v(B^{3})(1-p^2)^{{\frac{15}{2}}}\beta_{3} ^{5}dpm(d\omega)\\
 & = 3!v(B^{3})4\beta_4\beta_{3} ^{5}\int_0^1(1-p^2)^{{\frac{15}{2}}}dp
\end{align*}

\noindent We now solve for $v(B^{3})$. 

\noindent Using Lemma 1: $\int_0^1(1-p^2)^{{\frac{15}{2}}}dp = \dfrac{\beta_{16}}{2\beta_{15}}$, so
\begin{align*}
\beta_4^4 & = 3!v(B^{3})4\beta_4\beta_{3} ^{5}\dfrac{\beta_{16}}{2\beta_{15}}\\
2\beta_{15}\beta_4^4 & = 3!v(B^{3})4\beta_4\beta_{3} ^{5}\beta_{16}\\
v(B^{3}) & = \dfrac{2\beta_{15}\beta_4^4}{3!4\beta_4\beta_{3} ^{5}\beta_{16}}\\
 & = \dfrac{2\beta_{15}\beta_4^{3}}{4!\beta_{3}^{5}\beta_{16}}
\end{align*}
Using Lemma 2: $\beta_n = \dfrac{\pi^{\frac{n}{2}}}{(\frac{n}{2})!}$

\noindent Substituting Lemma 2 into the above, we get
\begin{align*}
v(B^{3}) & = \dfrac{2\dfrac{\pi^{\frac{15}{2}}}{\frac{15}{2}!}\dfrac{(\pi^{\frac{4}{2}})^{3}}{(\frac{4}{2}!)^{3}}}{4!\dfrac{(\pi^{\frac{3}{2}})^{5}}{(\frac{3}{2}!)^{5}}\dfrac{\pi^{\frac{16}{2}}}{\frac{16}{2}!}}\\\\
& = \dfrac{2\pi^{\frac{15}{2}}\pi^{\frac{12}{2}}(\frac{3}{2}!)^{5}\frac{16}{2}!}{ \frac{15}{2}!(\frac{4}{2}!)^{3}4!\pi^{\frac{15}{2}}\pi^{\frac{16}{2}}}
\end{align*}
We simplify to get
$$
v(B^{3}) = \dfrac{2\pi^{\frac{-4}{2}}(\frac{3}{2}!)^4\frac{3}{2}!\frac{16}{2}!\frac{4}{2}!}{ (\frac{4}{2}!)^4\frac{15}{2}!4!}
$$
Using the duplication formula $(m-\frac{1}{2})!=\dfrac{\pi^{\frac{1}{2}}(2m)!}{2^{2m}m!}$, we get
\begin{align*}
v(B^{3}) & = 2\pi^{\frac{-4}{2}}\dfrac{\dfrac{\pi^{\frac{4}{2}}(4!)^4}{2^{16}(\frac{4}{2}!)^4}\dfrac{\pi^{\frac{1}{2}}4!}{2^{4}\frac{4}{2}!}\frac{16}{2}!\frac{4}{2}!}{(\frac{4}{2}!)^4\dfrac{\pi^{\frac{1}{2}}16!}{2^{16}\frac{16}{2}!}4!}\\
& = \dfrac{2\pi^{\frac{-4}{2}}\pi^{\frac{4}{2}}(4!)^4\pi^{\frac{1}{2}}4!\frac{16}{2}!\frac{4}{2}!2^{16}\frac{16}{2}!}{2^{16}(\frac{4}{2}!)^42^{4}\frac{4}{2}!(\frac{4}{2}!)^4\pi^{\frac{1}{2}}16!4!}\\
& = \dfrac{(4!)^4\frac{16}{2}!\frac{16}{2}!}{(\frac{4}{2}!)^4 2^{3}(\frac{4}{2}!)^4 16!}\\
& = (\dfrac{4!}{\frac{4}{2}!\frac{4}{2}!})^4\dfrac{\frac{16}{2}!\frac{16}{2}!}{16!}2^{-3}\\
& = {\binom{4}{\frac{4}{2}}}^4{\binom{16}{\frac{16}{2}}}^{-1}2^{-3}.\\
\end{align*}

\begin{center}
\textbf{THE PROBLEM OF SYLVESTER FOR A UNIT BALL}
\end{center}

\noindent \textbf{THE PROBLEM:} Consider five points independently
and uniformly distributed over the volume of a convex region $K$ in
three-dimensional space. What is the probability that one of the
points falls within the tetrahedron formed by the other four points?\\

\noindent \textbf{THE SOLUTION IN THREE DIMENSIONS:}

The only way the region can fail to be convex is for one of the points to lie inside the tetrahedron formed by the other four points. There are five such configurations, depending on which point lies inside the tetrahedron, and they are mutually exclusive, because two points cannot both lie inside the tetrahedron at the same time. We may ignore the case where more than two of the points lie in a straight line, as this occurs with probability zero. The probability that a point lies inside the tetrahedron is given by $v(B^3)$.
\begin{align*}
v(B^{3}) & = {\binom{4}{\frac{4}{2}}}^4{\binom{16}{\frac{16}{2}}}^{-1}2^{-3}\\\\
& = {\binom{4}{2}}^4{\binom{16}{8}}^{-1} \frac{1}{8}\\\\
& = \left(\frac{4!}{2! 2!}\right)^4 \left(\frac{8! 8!}{16!}\right)\frac{1}{8}\\\\
& = \left(\frac{6^4}{8}\right)\left(\frac{1}{12870}\right) = \frac{9}{715}.
\end{align*}
\textbf{Answer to Sylvester Problem for a unit ball =} $1-5\left(\frac{9}{715}\right)= \mathbf{1-\dfrac{9}{143}} = \dfrac{134}{143}$\\

\subsection{Miscellaneous}

\noindent Here are the proofs of the lemmas and formulas we used above.
\begin{lem} \textbf{1.}
$\int_0^1(1-p^2)^{\frac{N-1}{2}}dp = \dfrac{\beta_N}{2\beta_N-1}$
\end{lem}
\begin{proof} Let $\beta_N$ be the volume of the unit ball in $\mathbb{R}^N$
\begin{align*}
\beta_N & = \underset{x_1^2+\cdots+x_N^2\leq1}{\int \cdots \int}dx_1dx_2 \cdots dx_N\\
& = \int_{-1}^1 dx \underset{x_2^2+\cdots+x_N^2\leq1-x^2}{\int \cdots \int} dx_2dx_2 \cdots dx_N
\end{align*}
The inner integral above is the volume of the ball of radius $(1-x^2)^{1/2}$ in $\mathbb{R}^{N-1}$. Hence we get
\begin{align*}
\beta_N & = \int_{-1}^1 [(1-x^2)^{\frac{1}{2}}]^{N-1}dx \beta_{N-1}\\
& = 2\int_0^1 (1-x^2)^{\frac{1}{2}{(N-1)}}dx \beta_{N-1}
\end{align*}
$$
\int_0^1 (1-x^2)^{\frac{1}{2}{(N-1)}}dx = \dfrac{\beta_N}{2\beta_{N-1}}.
$$
Let $p=x$, then
$$
\int_0^1 (1-p^2)^{\frac{N-1}{2}}dp = \dfrac{\beta_N}{2\beta_{N-1}}.
$$
\end{proof}

\begin{lem} \textbf{2.}
$\beta_n = \dfrac{\pi^{\frac{n}{2}}}{\frac{n}{2}!}$ for $n\geq1$.
\end{lem}
\begin{proof} A ball $B(x_0,r)$ of radius $r$ in $\mathbb{R}^n$ has volume $\beta_n r^n$ where $\beta_n$ is the volume of the unit ball.
 
Let $\Gamma(s)=\int_0^{\infty}y^{s-1}e^{-y}dy$, this is called the gamma function. Integrating by parts yields $\Gamma(s)= (s-1)\Gamma(s-1)$. From this, it is not hard to prove that $\Gamma(n)= (n-1)!$, $n \in \mathbb{N}$.  Because of this formula, the gamma function is often used to define $r!$ when $r$ is not an integer, i.e. $r!=\Gamma(r+1)$. 

Write the points in $\mathbb{R}^{n+1}$ as $(x,y)$ with $x\in \mathbb{R}^n$ and $y\in \mathbb{R}$, and let $D$ be the subset of $\mathbb{R}^{n+1}$ given by $D=\{(x,y):|x|^2 \leq y\}$. 

Integrate $e^{-y}$ over $D$, and use Fubini's theorem. In particular, Fubini's theorem states that the integral of $e^{-y}$ on $\mathbb{R}^{n+1}$ can be computed by iterating lower-dimensional integrals, and that the iterations can be taken in any order.
Integrating first over $x\in \mathbb{R}^n$, we get
$$
\int_{y=0}^{\infty} \int_{|x|^2 \leq y} e^{-y} dxdy = \int_0^{\infty} \beta_n y^{\frac{n}{2}}e^{-y}dy = \beta_n \Gamma(\frac{n}{2}+1)
$$
since the inner integral is over the ball $B(0,y^{\frac{1}{2}})$ of volume $\beta_n y^{\frac{n}{2}}$. On the other hand, integrating first over $y$, we get
$$
\int_{x\in \mathbb{R}^n} \int_{y=|x|^2}^{\infty} e^{-y} dydx = \int_{x\in \mathbb{R}^n} e^{-|x|^2} dx = \prod_{i=1}^n \int_{-\infty}^{\infty} e^{-x_i^2}dx_i = \left(\int_{-\infty}^{\infty} e^{-x^2}dx\right)^n.
$$
Fubini's theorem shows that the values above are equal, so we have
$$
\beta_n \Gamma\left(\frac{n}{2}+1\right)=\left(\int_{-\infty}^{\infty} e^{-x^2}dx\right)^n
$$
In a special case $n=2$, we have $\beta_2=\pi$ (area of an unit circle) and $\Gamma(2)=1$. Hence $\pi =(\int_{-\infty}^{\infty} e^{-x^2}dx)^2$, therefore $\int_{-\infty}^{\infty} e^{-x^2}dx=\pi^{\frac{1}{2}}$. So
$$
\beta_n \Gamma\left(\frac{n}{2}+1\right)=\pi^{\frac{n}{2}},
$$
i.e. 
$$
\beta_n =\dfrac{\pi^{\frac{n}{2}}}{\Gamma(\frac{n}{2}+1)}=\dfrac{\pi^{\frac{n}{2}}}{\frac{n}{2}!}.
$$
\end{proof}

To compute $\left(\dfrac{n}{2}\right)!$ for odd $n$, we have
\begin{lem} \textbf{3.} The duplication formula
$(m-\frac{1}{2})!=\dfrac{\pi^{\frac{1}{2}}(2m)!}{2^{2m}m!}$ for $m \in \mathbb{Z}$, $m \geq 0$.
\end{lem}
\begin{proof}
First we will show that $(-\frac{1}{2})! = \sqrt{\pi}$. 
 
Using the fact that $\Gamma(n)=\int_0^{\infty}y^{n-1}e^{-y}dy = (n-1)!$, we get that
$$
\left(-\frac{1}{2}\right)! = \left(\frac{1}{2}-1\right)! = \Gamma\left(\frac{1}{2}\right)
$$
$$
\Gamma\left(\frac{1}{2}\right)=\int_0^{\infty}y^{-\frac{1}{2}}e^{-y}dy
$$
Let $y=u^2$, so $dy = 2udu$
\begin{align*}
\Gamma\left(\frac{1}{2}\right) & = \int_0^{\infty}\frac{1}{u}e^{-u^2}2udu\\
& = 2\int_0^{\infty}e^{-u^2}du
\end{align*}
Let $I=\int_0^{\infty}e^{-u^2}du$. 
So $I^2 = \int_0^{\infty}e^{-x^2}dx \int_0^{\infty}e^{-y^2}dy$

\noindent By Fubini's Theorem, we get
$$
I^2 = \int_0^{\infty} \int_0^{\infty}e^{-(x^2+y^2)} dxdy
$$
We shall now change to polar coordinates to get
\begin{align*}
I^2 & = \int_0^{\frac{\pi}{2}} \int_0^{\infty}re^{-r} dr d\theta\\
& = \int_0^{\frac{\pi}{2}} \left[\frac{-e^{-u}}{2}\right]_0^{\infty} d\theta\\
& =\int_0^{\frac{\pi}{2}} \frac{1}{2} d\theta\\
& =\frac{\pi}{4}
\end{align*}
Hence $I = \frac{\sqrt{\pi}}{2}$
$$
\Gamma\left(\frac{1}{2}\right) = 2I = 2\frac{\sqrt{\pi}}{2} = \sqrt{\pi}.
$$
Now we will prove the duplication formula using induction
Let $m=0$, then $(-\frac{1}{2})! = \sqrt{\pi}$
Assume $(k-\frac{1}{2})!=\dfrac{\pi^{\frac{1}{2}}(2k)!}{2^{2k}k!}$ is true, then
\begin{align*}
\left(k+1-\dfrac{1}{2}\right)! & = \left(k+\dfrac{1}{2}\right)!\\
& = \left(k+\dfrac{1}{2}\right)\left(k+\dfrac{1}{2}-1\right)!\\
& = \left(k+\dfrac{1}{2}\right)\left(k-\dfrac{1}{2}\right)!\\
& = \left(k+\dfrac{1}{2}\right)\left(\dfrac{\pi^{1/2}(2k)!}{2^{2k}k!}\right)\\
& = \left(\dfrac{2k+1}{2}\right)\left(\dfrac{\pi^{1/2}(2k)!}{2^{2k}k!}\right)\left(\dfrac{2k+2}{2k+2}\right)\\
& = \dfrac{\pi^{1/2}(2k+2)!}{2^{2k+2}(k+1)!}\\
& = \dfrac{\pi^{1/2}(2(k+1))!}{2^{2(k+1)}(k+1)!}
\end{align*}
\end{proof}

Finally, note that if $n+1$ is odd, then
$(\frac{n+1}{2})!=(\frac{n+2}{2}-\frac{1}{2})!=\dfrac{\pi^{\frac{1}{2}}(n+2)!}{2^{n+2}(\frac{n+2}{2})!}$.
Thus
\begin{align*}
\binom{n+1}{\frac{1}{2}(n+1)} & = \dfrac{(n+1)!}{(\frac{n+1}{2})!(\frac{n+1}{2})!}\\
& = \dfrac{(n+1)!}{\dfrac{\pi^{\frac{1}{2}}(n+2)!}{2^{n+2}(\frac{n+2}{2})!}\dfrac{\pi^{\frac{1}{2}}(n+2)!}{2^{n+2}(\frac{n+2}{2})!}}\\
& =\dfrac{2^{2n+4}(\frac{n+2}{2})!(\frac{n+2}{2})!(n+1)!}{\pi(n+2)!(n+2)!}
\end{align*}
Hence we have that
$$
\binom{n+1}{\frac{1}{2}(n+1)} =\dfrac{2^{2n+4}}{\pi(n+2)}{\binom{n+2}{\frac{1}{2}(n+2)}}^{-1}.
$$

\section{Conclusions}

\noindent This thesis was initiated by trying to solve a Putnam problem which is a special case of Sylvester's four point problem. This problem was solved three different ways. The first found the expected area of the triangle formed by three of the four points in the three possible cases. In the second, we divided the interior of the circle into seven regions with three chords, then compared the expected area of the regions to find the expected area of the triangle. The third found the expected area of the triangle by finding the expected area of a triangle formed by two random points in a circle and a fixed point on the boundary.

In chapter 1, we proved Crofton's formula and Second Theorem, which is Sylvester's four point problem in dimension one.  Crofton's formula is the expected distance between two points chosen at random in a convex body. In this chapter, we developed Crofton's Second Theorem by describing the two dimensional results of Kingman's paper. $v(K)$ is the expected volume of the simplex where a point can be chosen inside that simplex, divided by the $n-$dimensional volume of $K$. Let $B^{n-1}$ be the unit ball in the $n-1$-dimension. We found $v(B)$ to be $\frac{1}{3}$.

In this thesis, we examine ``Random Secants of a Convex Body" in three and four dimensions to find $v(B^2)$ to be $\frac{35}{48\pi^2}$ and $v(B^3)$ to be $\frac{9}{715}$. With these results, we found that the solution to Sylvester's
four point problem in two dimensions is $1-\frac{35}{12\pi^2}$ and in three dimensions is $\frac{134}{143}$. Kingman generalized the results in his paper for a unit ball to be $v(B^{n-1}) = {\binom{n}{\frac{n}{2}}}^n{\binom{n^2}{\frac{n^2}{2}}}^{-1} 2^{-n+1}$

The Crofton formula has numerous applications. We will mention two more now: the Radon transform and the Busemann-Petty problem.

The Radon transform in two dimensions is the integral transform consisting of the integral of a function over the set of all lines. The Radon transform is an integral transform whose inverse is used to reconstruct images from medical CT scans. A technique for using Radon transforms to reconstruct a map of a planet's polar regions using a space craft in a polar orbit has also been devised.

The Busemann-Petty problem states that
\begin{quote}
\textit{``Let $K$ and $L$ be two origin-symmetric convex bodies in $\mathbb{R}^n$ such that the volume of $K$ intersected with every hyperplane through the origin is less than or equal to the corresponding volume of $L$ intersected with the same hyperplane. Does it follow that the volume of $K$ is less than or equal to the volume of $L$?"}
\end{quote}
The final answer to the problem is that the answer is ``yes" for $n \leq 4$ and ``no" for all larger values of $n$.

\end{document}